\documentclass[11pt]{article}
\newtheorem{theorem}{Theorem}
\newtheorem{lemma}{Lemma}
\newtheorem{proposition}{Proposition}
\newtheorem{corollary}{Corollary}
\newtheorem{definition}{Definition}
\newtheorem{remark}{Remark:}

\newcommand{\lab}[1]{\label{#1}}
\newcommand{\labs}[1]{\label{#1}}
\newcommand{\labe}[1]{\label{#1}}
\newcommand{\dating}[3]{\date{{\Large {\bf #1}}\vspace{5mm}\\{\Large#3}}\label{#2}}

\usepackage{amssymb}
\def\Bbb{\mathbb}

\newcommand{\Section}[1]{\section{#1}\setcounter{equation}{0}}

\def\Var{\mathop{\rm Var}\nolimits}

\newcommand{\cqfd}{\qquad\framebox[2.7mm]{\rule{0mm}{.7mm}}}
\newcommand{\bm}[1]{\mbox{\boldmath $#1$}}
\newcommand{\scal}[2]{\langle #1,#2\rangle}

\newcommand{\str}[1]{\rule{0mm}{#1mm}}
\newcommand{\st}{\strut}
\newcommand{\1}{1\hskip-2.6pt{\rm l}}

\begin{document}

\title{
\mbox{}\vspace{-10mm}\\{\Huge{\bf Model selection\vspace{2mm}
for density estimation with $\Bbb{L}_2$-loss}}\vspace{5mm}}
\author{
{\Large Lucien Birg\'{e}\vspace{3mm}}\\{\Large Universit\'e Paris VI\vspace{1mm}}\\
{\Large Laboratoire de Probabilit\'es et Mod\`eles Al\'eatoires\vspace{1mm}}\\
{\Large U.M.R.\ C.N.R.S.\ 7599}
}
\dating{\mbox{}}{April 2004}{19/01/2013}     

\maketitle
\footnotetext{

{\it \hspace{2.3mm}Key words and phrases.} Density estimation, $\Bbb{L}_2$-loss,
model selection, estimator selection, histograms.}

\begin{abstract}
We consider here estimation of an unknown probability density $s$ belonging to
$\Bbb{L}_2(\mu)$ where $\mu$ is a probability measure. We have at hand $n$ i.i.d.\
observations with density $s$ and use the squared $\Bbb{L}_2$-norm as our loss function.
The purpose of this paper is to provide an abstract but completely general method for
estimating $s$ by model selection, allowing to handle arbitrary families of finite-dimensional
(possibly non-linear) models and any $s\in\Bbb{L}_2(\mu)$. We shall, in particular,
consider the cases of unbounded densities and bounded densities with unknown
$\Bbb{L}_\infty$-norm and investigate how the $\Bbb{L}_\infty$-norm of $s$ may influence
the risk. We shall also provide applications to adaptive estimation and aggregation of 
preliminary estimators. One major technical tool of our approach is a proof of the existence of 
suitable tests between $\Bbb{L}_2$-balls with centers belonging to $\Bbb{L}_\infty$. Although of a purely theoretical nature, our method leads to results that cannot presently be reached by more concrete ones.
\end{abstract}

\Section{Introduction\labs{I}}
In this paper we shall deal with the problem of estimating an unknown density $s$ with respect to the measure $\mu$ on the measurable  space $({\cal X}, {\cal W})$ from an i.i.d.\ sample 
$\bm{X}=(X_1,\ldots,X_n)$ of random variables $X_i\in{\cal X}$ with distribution 
$P_s=s\cdot\mu$. We shall measure the quality of an estimator $\widehat{s}(X_1,\ldots,X_n)$ of $s$ by its quadratic risk $\Bbb{E}_s\left[d^2\left(\widehat{s},s\right)\right]$ for a suitable distance $d$, where $\Bbb{E}_s$ denotes the expectation when $s$ obtains. We shall denote by 
$\|\cdot\|_q$ the norm in $\Bbb{L}_q(\mu)$, omitting the subscript when $q=2$ for simplicity and by $d_2$ the distance in $\Bbb{L}_2(\mu)$: $d_2(t,u)=\|t-u\|$. For $1\le q\le+\infty$ we  consider the set $\overline{\Bbb{L}}_q$ of those densities with respect to $\mu$ that belong to $\Bbb{L}_q(\mu)$, i.e.
\[
\overline{\Bbb{L}}_q=\left\{t\in\Bbb{L}_q(\mu)\,\left|\,t\ge0\;\;\mbox{and}\;\;
\int t\,d\mu=1\right.\right\}.
\]
We shall also make use of the Hellinger distance $h$ and the variation distance $v$ given by
\[
h^2(t,u)=\frac{1}{2}\int\left(\sqrt{t}-\sqrt{u}\right)^2\,d\mu\qquad\mbox{and}\qquad
v(t,u)=\frac{1}{2}\int |t-u|\,d\mu.
\]

When $s$ is assumed to belong to the metric space $(M,d)$, a common method of estimation that can be called model-based estimation chooses a subset $\overline{S}$ of $M$ and an estimation method which results in an estimator that automatically belongs to $\overline{S}$. Of this type is the maximum likelihood estimator over $\overline{S}$, for instance. When the distance $d$ is either $h$ or $v$, it is possible to get very general risk bounds for some special estimators based on finite-dimensional models. Indeed, if we choose for $\overline{S}$ a model with a metric dimension (to be defined in Section~\ref{I8} below and generalizing to subsets of metric spaces the usual dimension of a finite-dimensional linear space) bounded by $D$, one can design an estimator $\widetilde{s}$ with values in $\overline{S}$ satisfying, whatever the true unknown density $s$,
\begin{equation}
\Bbb{E}_s\left[d^2\left(\widetilde{s},s\right)\right]\le 
C\left[\inf_{t\in\overline{S}}d^2(s,t)+n^{-1}D\right],
\labe{Eq-M15a}
\end{equation}
where $C$ denotes a universal constant (independent of $n$, $s$ and $\overline{S}$). When 
$s$ belongs to $\overline{S}$, (\ref{Eq-M15a}) provides the following upper bound for the minimax risk over $\overline{S}$:
\begin{equation}
\inf_{\widehat{s}}\sup_{s\in\overline{S}}\Bbb{E}_s\left[d^2\left(\widehat{s},s\right)\right]\le 
Cn^{-1}D,
\labe{Eq-M15b}
\end{equation}
where the infimum is over all possible estimators $\widehat{s}$, a result which actually dates back to Le Cam (1973) for the Hellinger distance.

Nevertheless, the square of the $\Bbb{L}_2$-distance $d_2$ has been much more popular in the past, as a loss function for density estimation, than either the Hellinger or variation distances, mainly because of its simplicity due to the classical ``squared bias plus variance" decomposition of the risk. But, although hundreds of papers have been devoted to the derivation of risk bounds for various specific estimators, we do not know of any {\em universal} bound for the risk similar to 
(\ref{Eq-M15a}) when $d=d_2$, universal meaning here valid for any model $\overline{S}$ with a metric dimension bounded by $D$ and any density $s\in\Bbb{L}_2(\mu)$, only partial results valid for some special cases being available. This is actually not surprising for the following reason. While $h$ and $v$ are distances between probabilities so that $h(s,t)=h(P_s,P_t)$ is independent of the choice of the underlying dominating measure $\mu$, this is definitely not the case of the $\Bbb{L}_2$-distance between densities which depends on the choice of $\mu$ and is not a distance between probabilities. Given a probability $P$ and a dominating measure $\mu$, even the fact that $dP/d\mu$ belong or not to $\Bbb{L}_2(\mu)$ depends on $\mu$. Further remarks on this subject can be found in Devroye and Gy\"orfi (1985) and Devroye (1987).

It is indeed the distortion between the Hellinger and $\Bbb{L}_2$-distances that explains the problems that may occur when we use the $\Bbb{L}_2$-risk as can be shown by the following elementary computations. When $t$ and $u$ are bounded by $L$,
\begin{equation}
\|t-u\|^2=\int\left(\sqrt{t}-\sqrt{u}\right)^2\left(\sqrt{t}+\sqrt{u}\right)^2d\mu\le4L
\int\left(\sqrt{t}-\sqrt{u}\right)^2d\mu=8Lh^2(t,u).
\labe{Eq-HL}
\end{equation}
Although this is only an upper bound, there are situations where it is rather sharp 
as in the following case. Let $\mu$ be the Lebesgue measure on $[0,1]$, 
$t=L\1_{[0,a)}+(1-aL)(1-a)^{-1}\1_{[a,1]}$ with $L>1$ and $0<a<L^{-1}$, and $u(x)=t(1-x)$ for $0\le x\le1$. Then $\|t\|_\infty=\|u\|_\infty=L$ and
\[
\|t-u\|^2=2a\left(L-\frac{1-aL}{1-a}\right)^2=2a\frac{(L-1)^2}{(1-a)^2},
\]
while
\[
h^2(t,u)=a\left(\sqrt{L}-\sqrt{\frac{1-aL}{1-a}}\right)^2=\frac{a}{1-a}\left(\frac{L-1}{\sqrt{L(1-a)}
+\sqrt{1-aL}}\right)^2\le\frac{a(L-1)^2}{L(1-a)^2}.
\]
Therefore $\|t-u\|^2\ge2Lh^2(t,u)$. If $a$ is chosen in such a way that $h^2(t,u)=(4n)^{-1}$, it follows from Le Cam (1973) (see Proposition~\ref{P-Mino} below) that one cannot test between $t$ and $u$ with $n$ i.i.d.\ observations and small errors, i.e.\ one cannot distinguish between $t$ and $u$ with only $n$ observations. As a consequence the minimax risk over the set $\{t,u\}$ will be of order $n^{-1}$ when the loss function is the squared Hellinger distance while it will be of order $L n^{-1}$ if the loss function is the squared 
$\Bbb{L}_2$-distance. 

\subsection{The example of projection estimators\labs{I2}}

\subsubsection{The special case of histograms\labs{I2a}}
A simple illustration of the difference that occurs when one computes risk bounds using the $\Bbb{L}_2$-distance rather than the Hellinger distance is provided by the case of histograms.
Assuming that $\mu$ is a finite measure and given a finite partition ${\cal I}=\{I_1,\ldots,I_k\}$ of ${\cal X}$ with $\mu(I_j)=l_j>0$ for $1\le j\le k$, the histogram $\widehat{s}_{{\cal I}}$ based on this partition is defined by
\begin{equation}
\widehat{s}_{{\cal I}}(X_1,\ldots, X_n)=\frac{1}{nl_j}\sum_{j=1}^kN_j\1_{I_j},
\quad\mbox{with } N_j=\sum_{i=1}^n\1_{I_j}(X_i).
\labe{Eq-Hist}
\end{equation}
Let 
\[
p_j=\int_{I_j}s\,d\mu,\quad\overline{s}_{{\cal I}}=\sum_{j=1}^k\frac{p_j}{l_j}\1_{I_j},\qquad
\overline{S}_{{\cal I}}=
\left\{\left.\sum_{j=1}^k\beta_j\1_{I_j}\,\right|\beta_j\in\Bbb{R}\mbox{ for }1\le j\le k\right\}
\]
and $\overline{S}_{{\cal I}}^0$ be the convex set $\overline{S}_{{\cal I}}\cap\overline{\Bbb{L}}_1$. If $s\in\Bbb{L}_2(\mu)$, then $\overline{s}_{{\cal I}}\in\overline{S}_{{\cal I}}^0$ is the orthogonal projection of $s$
onto the $k$-dimensional linear space $\overline{S}_{{\cal I}}$ spanned by the functions 
$\1_{I_j}$ and onto $\overline{S}_{{\cal I}}^0$ as well. It follows that $\widehat{s}_{{\cal I}}$ is an estimator based on the model $\overline{S}_{{\cal I}}^0$.

Choosing $d_2^2$ as our loss function, we derive that
\begin{equation}
\Bbb{E}_s\left[\|\widehat{s}_{{\cal I}}-s\|^2\right]=\|\overline{s}_{{\cal I}}-s\|^2+
\frac{1}{n}\sum_{j=1}^k\frac{p_j(1-p_j)}{l_j}.
\labe{Eq-his1}
\end{equation}
The simplest and most common situation is the one of {\em regular} histograms for which all $l_j$ are equal to
$k^{-1}$. In this case we derive from (\ref{Eq-his1}) and a convexity argument that
\begin{equation}
\Bbb{E}_s\left[\|\widehat{s}_{{\cal I}}-s\|^2\right]\le\|\overline{s}_{{\cal I}}-s\|^2+n^{-1}(k-1).
\labe{Eq-his2}
\end{equation}
We then get a risk bound which is the sum of the square of the bias and $n^{-1}$ times the dimension of the model, i.e.\ the number of parameters ($p_1,\ldots,p_{k-1}$) which are needed to describe an element of the model $\overline{S}_{{\cal I}}^0$. It can therefore be viewed as an analogue of (\ref{Eq-M15a}).

In the general situation of unequal values of the $l_j$ the previous elementary argument does not work but, if $s\in\Bbb{L}_\infty(\mu)$ with norm $\|s\|_\infty$, the quadratic risk of $\widehat{s}_{{\cal I}}$ can alternatively be bounded by
\begin{equation}
\Bbb{E}_s\left[\|\widehat{s}_{{\cal I}}-s\|^2\right]\le
\|\overline{s}_{{\cal I}}-s\|^2+n^{-1}(k-1)\|s\|_\infty,
\labe{Eq-his5}
\end{equation}
which is much worse than (\ref{Eq-his2}) when $\|s\|_\infty$ is large. This bound may be far from sharp for a given $s$ but it is essentially unimprovable if we want it to hold for arbitrary partitions ${\cal I}$ with $k$ elements and any $s\in\overline{\Bbb{L}}_\infty$ as shown by the following example. Define the partition ${\cal I}$  on ${\cal X}=[0,1]$ by $I_j=[(j-1)\alpha,j\alpha)$ for $1\le j<k$ and 
$I_k=[(k-1)\alpha,1]$ with $0<\alpha<(k-1)^{-1}$. Set $s=s_{{\cal I}}=
[(k-1)\alpha]^{-1}\left[1-\1_{I_k}\right]$. Then $p_j=(k-1)^{-1}$ for  $1\le j<k$, $s=\overline{s}_{{\cal I}}$ 
(a case of no bias) and it follows from (\ref{Eq-his1}) that
\begin{equation}
\Bbb{E}_s\left[\|\widehat{s}_{{\cal I}}-s\|^2\right]=\frac{k-2}{(k-1)\alpha n}
=\frac{(k-2)\|s\|_\infty}{n}.
\labe{Eq-his4}
\end{equation}
This shows that there is little space for improvement in (\ref{Eq-his5}) and that there are cases when the quadratic risk based on $d_2$ does involve the $\Bbb{L}_\infty$-norm of $s$. It also demonstrates that there is no hope to bound the risk of an histogram $\widehat{s}_{{\cal I}}$ based on an arbitrary  partition ${\cal I}$ with cardinality $|{\cal I}|=k$ by an analogue of (\ref{Eq-his2}). Indeed, letting 
$\alpha$ go to zero in (\ref{Eq-his4}) shows that
\[
\sup_{\{{\cal I}\,|\,|{\cal I}|=k\}}\sup_{s\in\overline{S}_{{\cal I}}^0}
\Bbb{E}_s\left[\|\widehat{s}_{{\cal I}}-s\|^2\right]=+\infty.
\]

If, instead, we use as our loss function the squared Hellinger distance $h$ as we previously did we get an analogue of (\ref{Eq-his2}) and (\ref{Eq-M15a}), namely
\[
\Bbb{E}_s\left[h^2(\widehat{s}_{\cal I},s)\right]\le h^2(s,\overline{s}_{{\cal I}})+(k-1)/(2n),
\]
whatever the partition ${\cal I}$ of cardinality $k$ and the density $s$, as shown in Birg\'e and Rozenholc (2006). A similar result holds if $d=v$. In both cases, whatever the partition ${\cal I}$, we can bound the risk by a universal constant times the sum of the squared of the bias and $n^{-1}$ times the size of the partition. This is a bound of the form (\ref{Eq-M15a}), since $|{\cal I}|$ is the dimension of our model, i.e.\ the linear space generated by the functions $\1_{I_j}$ to which $\widehat{s}_{\cal I}$ belongs.

\subsubsection{Projection estimators\labs{I2b}}
More generally, instead of the linear space generated by the functions $\1_{I_j}$, $1\le j\le k$, we can take as a model for estimating $s$ any $k$-dimensional linear subspace $\overline{S}$ of 
$\Bbb{L}_2(\mu)$. Given an orthonormal basis $(\varphi_1,\ldots,\varphi_k)$ of 
$\overline{S}$ the projection $\overline{s}$ of $s$ onto $\overline{S}$ can be written
$\overline{s}=\sum_{j=1}^k\beta_j\varphi_j$. The estimation method of Cencov (1962) consists in replacing each coefficient $\beta_j=\int\varphi_js\,d\mu$ in this expansion by its empirical version $\widehat{\beta}_j=n^{-1}\sum_{i=1}^n\varphi_j(X_i)$. This results in the so called {\em projection estimator} $\widehat{s}=\sum_{j=1}^k\widehat{\beta}_j\varphi_j$ (which in general is not a density) with risk
\begin{eqnarray}
\Bbb{E}_s\left[\|\widehat{s}-s\|^2\right]&=&\|\overline{s}-s\|^2+n^{-1}\sum_{j=1}^k
\Var_s\left(\st\varphi_j(X_1)\right)
\labe{Eq-proj0}\\
&\le&\|\overline{s}-s\|^2+n^{-1}\int\left[\sum_{j=1}^k\varphi_j^2(x)\right]s(x)\,d\mu(x)
\labe{Eq-proj1}\\
&\le&\|\overline{s}-s\|^2+kn^{-1}\min\left\{k^{-1}\left\|\sum_{j=1}^k\varphi_j^2\right\|_\infty;
\|s\|_\infty\right\}.
\labe{Eq-proj2}
\end{eqnarray}
The histogram based on the partition ${\cal I}$ is merely the projection estimator corresponding to 
$\varphi_j=l_j^{-1/2}\1_{I_j}$. For regular histograms, $l_j=k^{-1}$ and $\left\|\sum_{j=1}^k\varphi_j^2\right\|_\infty=k$.

It has been shown in Birg\'e and Massart (1998) that the quantity $\left\|\sum_{j=1}^k\varphi_j^2\right\|_\infty$ only depends on $\overline{S}$ and not on the choice of the basis. For an arbitrary subset of a uniformly bounded basis like the trigonometric basis, it is bounded by $Ck$ for a constant $C$ depending on the basis only. In such a case we get a risk bound of the form $\|\overline{s}-s\|^2+n^{-1}Ck$ which does not involve $\|s\|_\infty$. If we use the projection onto the first $k$ elements of a wavelet basis, the bound $Ck$ still holds but this is not true any more if we project onto an arbitrary subset with $k$ elements of the same wavelet basis. We end up with the same dichotomy we found between regular and irregular histograms: bounding the variance term of the risk is sometimes straightforward, leading to a bound of the form $n^{-1}Ck$ but for some other models the risk bound we derive involves $\|s\|_\infty$.

A similar difficulty occurs in more complex examples, for instance for the estimators that are considered by Reynaud-Bouret, Rivoirard and Tuleau-Malot (2011). Their Theorem~1 leads to a risk bound that also involves a variance term depending on the unknown density $s$ which is the analogue of $\sum_{j=1}^k\Var_s\left(\st\varphi_j(X_1)\right)$. In some cases, this term can be bounded independently of $s$ but in some other cases this bounding involves $\|s\|_\infty$.

It follows from these illustrations that it does not seem easy to get an analogue of (\ref{Eq-M15a}) in full generality when $d$ is the $\Bbb{L}_2$-distance (at first sight it may sometimes work and sometimes not, depending on the type of model we use). It will be the subject of our next section to formally prove that a general result of the form (\ref{Eq-M15a}) cannot hold when $d=d_2$.

\subsection{Model based estimation\labs{I8}}
As we already mentioned, a common method for estimating $s$ consists in choosing a particular subset $\overline{S}$ of $(M,d)$ that we shall call a {\em model} for $s$ and design an estimator with values in $\overline{S}$. Let us set $M=\overline{\Bbb{L}}_1$ and choose for $d$ either the Hellinger distance $h$ or the variation distance $v$. It follows from Le Cam (1973, 1975, 1986) and subsequent results by Birg\'e (1983, 2006a) that the risk of suitably designed estimators with values in $\overline{S}$ is the sum of two terms, an {\em approximation term} depending on the distance from $s$ to $\overline{S}$ and an {\em estimation term} depending on the metric dimension of the model $\overline{S}$ which can be defined as follows.

\begin{definition}\lab{D-Bmdim}
Let $\overline{S}$ be a subset of some metric space $(M,d)$ and let ${\cal B}_d(t,r)$ denote
the open ball of center $t$ and radius $r$ with respect to the metric $d$. Given $\eta>0$, a subset 
$S_\eta$ of  $M$ is called an $\eta$-net for $\overline{S}$ if, for each $t\in\overline{S}$, one can find $t'\in S_\eta$ with $d(t,t')\le\eta$.

We say that $\overline{S}$ has {\em a metric dimension bounded by} $D\ge0$ if, for every
$\eta>0$, there exists an $\eta$-net $S_\eta$ for $\overline{S}$ such that
\begin{equation}
|S_\eta\cap{\cal B}_d(t,x\eta)|\le\exp\left[Dx^2\right]\quad\mbox{for
all }x\ge2\mbox{ and }t\in M.
\labe{Eq-M49}
\end{equation}
\end{definition}
%
%
\begin{remark}
One can always assume that $S_\eta\subset\overline{S}$ at the price of replacing $D$ by
$25D/4$ according to Proposition~7 of Birg\'e (2006a).
\end{remark}
When $(M,d)$ is a normed linear space, typical examples of sets with metric dimension bounded by 
$D$ are subsets of $2D$-dimensional linear subspaces of $M$, as shown in Birg\'e (2006a), where the following generalization of Le Cam (1973) is also proven. 
%
\begin{proposition}\lab{P-minmax}
Assume that we observe $n$ i.i.d.\ random variables with unknown distribution $P_s$,
$s\in\left(\overline{\Bbb{L}}_1,d\right)$, $d$ being either the Hellinger distance $h$ or the variation distance $v$, and that we have at disposal a subset $\overline{S}$ of $\overline{\Bbb{L}}_1$  with metric dimension bounded by $D\ge1/2$. One can build an estimator $\widetilde{s}$ with values in 
$\overline{S}$ such that, for all $s\in\overline{\Bbb{L}}_1$ and some universal constant $C$,
\[
\Bbb{E}_s\left[d^2\left(\widetilde{s},s\right)\right]\le 
C\left[\inf_{t\in\overline{S}}d^2(s,t)+n^{-1}D\right]\quad\mbox{hence}\quad
\sup_{s\in\overline{S}}\Bbb{E}_s\left[d^2\left(\widetilde{s},s\right)\right]\le Cn^{-1}D.
\]
\end{proposition}
The risk bounds (\ref{Eq-M15a}) and (\ref{Eq-M15b}) that we mentioned earlier actually derive from this proposition.

\subsection{Some negative results for the $\Bbb{L}_2$-loss\labs{R9}}
Unfortunately, the analogue of Proposition~\ref{P-minmax} when we deal with arbitrary densities and models belonging to $\overline{\Bbb{L}}_2$ and set $d=d_2$ cannot be true. To see this, let us take ${\cal X}=[0,1]$, $\mu$ the Lebesgue measure on ${\cal X}$ and assume that, for some $\overline{S}$ with metric dimension bounded by $D\ge1/2$,
\[
\inf_{\widehat{s}}\sup_{s\in\overline{S}}\Bbb{E}_s\left[d_2^2\left(\widehat{s}(\bm{X}),s\right)\right]=cn^{-1}D
\]
for some $c>0$. For $\lambda>1$, consider the mapping $G_\lambda$ between elements of  
$\overline{\Bbb{L}}_2$ given by $G_\lambda(s)(x)=\lambda s(\lambda x)
\1_{[0,\lambda^{-1}]}(x)$. Then $d_2\left(G_\lambda(t)-G_\lambda(u)\right)=
\lambda^{1/2}d_2(t,u)$. This implies that $G_\lambda\left(\overline{S}\right)$ has the same metric dimension as $\overline{S}$, therefore bounded by $D$. Moreover, any estimator $\widehat{s}$ of $s$ can be turned into an estimator $G_\lambda\left(\widehat{s}\right)$ for $G_\lambda(s)$ and vice-versa, so that the minimax risk over $G_\lambda\left(\overline{S}\right)$ is $\lambda cn^{-1}D$. Since $\lambda$ can be arbitrary large, the bound (\ref{Eq-M15b}) cannot be universally true. 

The fact that the $\Bbb{L}_\infty$-norm of $s$ may come into the risk based on $\Bbb{L}_2$-loss, as we noticed when studying histograms on irregular partitions, is actually not due to the use of specific estimators like histograms but it is a more general phenomenon as shown by another negative result provided by Proposition~4 of Birg\'e (2006b) that we recall below for the sake of completeness.
%
\begin{proposition}\lab{P-Contrex1}
For each $L>0$ and each integer $D$ with $1\le D\le 3n$, one can find a finite set
$\overline{S}$ of densities with the following properties:

i) it is a subset of some $D$-dimensional affine subspace of $\Bbb{L}_2([0,1],dx)$ with a 
metric dimension bounded by $D/2$;

ii) $\sup_{s\in\overline{S}}\|s\|_\infty\le L+1$;

iii) for any estimator $\widehat{s}(X_1,\ldots,X_n)$ belonging to $\Bbb{L}_2([0,1],dx)$ and 
based on an i.i.d.\ sample with density $s\in\overline{S}$,
\begin{equation}
\sup_{s\in\overline{S}}\Bbb{E}_s\left[\|\widehat{s}-s\|^2\right]>0.0139DLn^{-1}.
\labe{Eq-Z2}
\end{equation}
\end{proposition}
It follows from this lower bound that the best {\em universal} risk bound one can expect to prove for an estimator $\widehat{s}$ with values in an arbitrary model $\overline{S}$ with metric dimension bounded by $D$ when $s$ is arbitrary in $\overline{\Bbb{L}}_\infty$ is
\begin{equation}
\Bbb{E}_s\left[d_2^2\left(\widehat{s},s\right)\right]\le C\left[\inf_{t\in\overline{S}}
d_2^2(s,t)+n^{-1}D\|s\|_\infty\right].
\labe{Eq-M15c}
\end{equation}
The situation becomes worse when $s\not\in\Bbb{L}_\infty(\mu)$ or if
$\sup_{s\in\overline{S}}\|s\|_\infty=+\infty$. It may even happen that, whatever the estimator 
$\widehat{s}$, $\sup_{s\in\overline{S}}\Bbb{E}_s\left[d^2\left(\widehat{s},s\right)\right]$ be infinite even if $\overline{S}\subset\Bbb{L}_2(\mu)$ has a bounded metric dimension as shown by the following lower bound, to be proved in Section~\ref{P0}.
%
\begin{proposition}\lab{P-Contrex2}
Let $\overline{S}=\{s_\theta,\, 0<\theta\le1/3\}$ be the set of densities with
respect to the Lebesgue measure on $[0,1]$ given by
\[
s_\theta=\theta^{-2}\1_{[0,\theta^3]}+\left(\theta^2+\theta+1\right)^{-1}
\1_{(\theta^3,1]}.
\]
If  we have at disposal $n$ i.i.d.\ observations with density $s_\theta\in\overline{S}$, we 
can build an estimator $\widetilde{s}_n$ such that  $\sup_{0<\theta\le1/3}
\Bbb{E}_{s_\theta}\left[nh^2(s_\theta,\widetilde{s}_n)\right]\le C$ for some $C$ independent
of $n$. On the other hand, although the metric dimension of $\overline{S}$ with respect to
the distance $d_2$ is bounded by $2$, $\sup_{0<\theta\le1/3}\Bbb{E}_{s_\theta}
\left[\|s_\theta-\widehat{s}_n\|^2\right]=+\infty$, whatever $n$ and the estimator $\widehat{s}_n$.
\end{proposition}
These three counter-examples show that there is absolutely no hope that 
Proposition~\ref{P-minmax}  could be true when $d=d_2$. They also suggest that it is impossible to build a general theory of model selection (or even of estimation based on one single model) with $\Bbb{L}_2$-loss without taking the $\Bbb{L}_\infty$-norm into account, even if there do exist some special situations, like for regular histograms, for which the introduction of the 
$\Bbb{L}_\infty$-norm is superfluous and Bound~(\ref{Eq-M15c}) can be substantially improved.

\subsection{About this paper\labs{I3}}
We have seen in the previous section that Proposition~\ref{P-minmax}, which deals with one single model $\overline{S}$, cannot be true when $d$ is the $\Bbb{L}_2$-distance and the situation is obviously worse for model selection among many models. The general results about model selection or estimator aggregation (that can be viewed as a special case of model selection as explained in Section~9 of Birg\'e (2006a)) which are valid when $d=h$ or $v$ and described precisely in Theorem~\ref{T-sevmodel} below cannot hold in full generality when $d=d_2$. Of course, they may hold in some specific situations or under some additional restrictions and many results have already been obtained in this direction but there is presently no general theory for model selection available when the risk is based on the $\Bbb{L}_2$-distance. This major difference between the $\Bbb{L}_2$-distance and 
$\Bbb{L}_1$ or Hellinger distance was the main motivation to write this paper.

Our purpose, in the remainder of this paper, will be to explain to what extent the general theory for model selection which has been developed in Birg\'e (2006a) for $d=h$ or $v$ can be rescued when $d=d_2$ with the additional introduction of $\Bbb{L}_\infty$-norms in the procedures, even when the density $s$ does not belong to $\Bbb{L}_\infty(\mu)$, and what type of results about adaptation in Besov spaces can be derived from this general approach.
In particular, for the case of a single model, this will lead to a generalized version of  
(\ref{Eq-M15b}) that can also handle the case of $s\not\in\Bbb{L}_\infty$. When 
$s\in\Bbb{L}_\infty$ (with an unknown value of $\|s\|_\infty$), the risk bounds we get completely parallel (apart from some constants depending on $\|s\|_\infty$) those obtained 
for estimating $s$ when $d=h$ or $v$. 

In order to achieve our goal, we shall use the construction of what we have called T-estimators in Birg\'e (2006a). These estimators are based on suitable tests between balls with respect to the relevant distance $d$. In the i.i.d.\ case when $d=h$, tests between Hellinger balls were constructed quite a long time ago by Le Cam and there are many other frameworks for which such tests exist; see Birg\'e (2012) for various examples. In order to apply our construction to the case of $d=d_2$ we shall have to derive suitable tests between 
$\Bbb{L}_2$-balls. 

In the next section we shall recall general results for model selection, based on what we have called T-estimators in Birg\'e (2006a), that hold when $d=h$ or $v$ and what is presently known when $d=d_2$. Section~\ref{U} will be devoted to the statement of the main theorems and we shall give a few applications of them, in particular to aggregation of preliminary estimators and estimation of densities belonging to Besov spaces, in Section~\ref{X}. We shall explain precisely the construction of our estimators, which can be viewed as a modification of T-estimators, in Section~\ref{T}.  The last section will be devoted to the most technical proofs.

\Section{Model selection\labs{R6}}
Let us now go back to histograms. As we noticed, the risk bound we get heavily depends on the choice of the partition. If we have at disposal a finite (although possibly very large) family 
$\{{\cal I}_m,m\in{\cal M}\}$ of finite partitions of ${\cal X}$ with respective cardinalities 
$|{\cal I}_m|$, we can consider the corresponding families of models $\{\overline{S}_{{\cal I}_m}, m\in{\cal M}\}$ and histogram estimators $\{\widehat{s}_{{\cal I}_m}, m\in{\cal M}\}$. It is then natural to try to find a partition in the family which leads, at least approximately, to the minimal risk $\inf_{m\in{\cal M}}\Bbb{E}_s\left[\|\widehat{s}_{{\cal I}_m}-s\|^2\right]$. But one cannot select
such a partition from either (\ref{Eq-his1}) or (\ref{Eq-his5}) since the risk depends on the unknown density $s$ via $\overline{s}_{{\cal I}_m}$. Methods of {\em model or estimator selection} base the choice of a suitable partition ${\cal I}_{\widehat{m}}$ with $\widehat{m}=
\widehat{m}(X_1,\ldots,X_n)$ on the observations. 

This problem of partition selection is actually a particular case of model selection. Going back to the general framework of Section~\ref{I8} with $d=h$ or $v$, we can consider a family of models 
$\left\{\overline{S}_m,m\in{\cal M}\right\}$, each one with metric dimension bounded by $D_m$ so that it leads to an estimator $\widehat{s}_m$ with a risk bounded, according to 
Proposition~\ref{P-minmax}, by $C\left[\inf_{t\in\overline{S}_m}d^2(s,t)+n^{-1}D_m\right]$. Since the bias term $\inf_{t\in\overline{S}_m}d^2(s,t)$ is unknown, it is impossible to decide which $m$ leads to the best bound and a natural problem is to design a method for choosing a value 
$\widehat{m}(X_1,\ldots,X_n)$ of $m$ from the observations in order to minimize the risk bound.
There is actually a solution to this problem which is provided by the following result from Birg\'e (2006a).

\begin{theorem}\lab{T-sevmodel} 
Let $\bm{X}=(X_1,\ldots,X_n)$ be an i.i.d.\ sample with unknown density $s$ belonging to
$\overline{\Bbb{L}}_1$, $d$ be either $h$ or $v$ and $\left\{\overline{S}_m,m\in{\cal M}\right\}$ a finite or countable family of subsets of $\overline{\Bbb{L}}_1$ with metric dimensions bounded by $D_m\ge1/2$ respectively. Let the nonegative weights $\Delta_m, m\in{\cal M}$ satisfy
\begin{equation}
\sum_{m\in{\cal M}}\exp[-\Delta_m]=\Sigma'<+\infty.
\labe{Eq-S2}
\end{equation}
Then there exists a universal constant $C$ and an estimator
$\widetilde{s}(\bm{X})$ such that, for any  $s\in\overline{\Bbb{L}}_1$,
\begin{equation}
\Bbb{E}_s\left[d^2\left(\widetilde{s},s\right)\right]\le C(1+\Sigma')\inf_{m\in{\cal M}}
\left[\inf_{t\in\overline{S}_m}d^2(s,t)+n^{-1}\max\left\{D_m;\Delta_m\right\}\right].
\labe{Eq-M15f}
\end{equation}
\end{theorem}
Proposition~\ref{P-minmax} simply follows by setting ${\cal M} =\{0\}$, $\overline{S}_0=\overline{S}$, $D_0=D$ and $\Delta_0=1/2$. One should notice here the Bayesian role of the weights $\Delta_m$. The choice of weights that satisfy (\ref{Eq-S2}) amounts to putting a prior positive measure with total mass $\Sigma'$ on ${\cal M}$ or equivalently on the collection of models with a measure $\exp[-\Delta_m]$ for the model $\overline{S}_m$.

\subsection{What is presently known\labs{R2}}
There exists a considerable amount of literature dealing with problems of model or
estimator selection. Most of it is actually devoted to the analysis of Gaussian problems,
or regression problems, or density estimation with either Hellinger or Kullback loss and it is
not our aim here to review this literature. Only a few papers are actually devoted to our
subject, namely model or estimator selection for estimating densities with
$\Bbb{L}_2$-loss, and we shall therefore concentrate on these papers only. They  can
roughly be divided into three groups: the ones dealing with penalized projection estimators, the ones that study aggregation by selection of preliminary estimators and the more specific ones which use methods based on the thresholding of empirical coefficients within a
given wavelet basis. The last ones, which are especially designed for the estimation of densities belonging to various kinds of Besov spaces, are typically not advertised as dealing with model selection but, as explained for instance in Section~5.1.2 of Birg\'e and Massart (2001), can be viewed as very special instances of model selection methods for models that are spanned by some finite subsets of an orthonormal wavelet basis. They are definitely not general methods of model selection (i.e.\ which handle arbitrary densities and families of models) but specific ones dealing only with some special families of models and targeted to estimate special densities.

All these papers have in common the fact that they require more or less severe restrictions on the families of models or densities to be estimated. For instance, aggregation of estimators by selection only deals with models which are singletons while thresholding of wavelet coefficients amounts to deal with models which are linear spaces spanned by finite subsets of a wavelet basis. Moreover, apart from a few special cases to be mentioned below, they typically assume that $s\in\Bbb{L}_\infty(\mu)$ with a known or estimated bound for 
$\|s\|_\infty$. 

In order to see how such methods apply to the problem of partition selection for histograms that we mentioned at the beginning of Section~\ref{R6}, let us be more specific and assume that 
${\cal X}=[0,1]$, $\mu$ is the Lebesgue measure and ${\cal N}=\{j/(N+1), 1\le j\le N\}$ for some (possibly very large) positive integer $N$. For any subset $m$ of ${\cal N}$, we denote by ${\cal I}_m$ the partition of ${\cal X}$ generated by the intervals with set of endpoints 
$m\cup\{0,1\}$ and we set $\overline{S}_m=\overline{S}_{{\cal I}_m}$ and 
$\widehat{s}_m=\widehat{s}_{{\cal I}_m}$. This leads to a family of partitions ${\cal M}$ with cardinality $2^N$ and to the corresponding families of linear models $\left\{\overline{S}_m, 
m\in{\cal M}\right\}$ and related histogram estimators $\{\widehat{s}_m, m\in{\cal M}\}$. Then all models $\overline{S}_m$ are subsets of the largest one $\overline{S}_{{\cal N}}$. Given a sample $X_1,\ldots,X_n$ with unknown density $s$, which partition ${\cal I}_{\widehat{m}}$ with $\widehat{m}=\widehat{m}(X_1,\ldots,X_n)$ depending on the observations
should we choose to estimate $s$ and what sort of risk bound could we derive for the resulting estimator~? 

Since subset selection within a given basis applies here only when $N=2^K-1$ and we use the Haar basis, we shall only consider this particular case in order to be able to deal with the three above-mentioned methods, keeping in mind that the third one does not apply to arbitrary values of $N$. In this case, $\overline{S}_{{\cal N}}$ is the linear span of the $2^K$ first elements of the Haar basis. To each non-empty subset $q$ of these $2^K$ elements, we can associate its linear span $\overline{S}'_q$ and the family of linear models $\left\{\overline{S}'_q, q\in{\cal Q}\right\}$ where ${\cal Q}$ denotes the set of those $q$. To each $\overline{S}'_q$ corresponds a projection estimator (as defined in Section~\ref{I2b}) $\widehat{s}_q$ which looks like a histogram estimator (piecewise constant) and one can consider the problem of selecting an optimal value $\widehat{q}$ of $q$, for instance by a proper thresholding of the empirical coefficients. One should nevertheless keep in mind that the two problems (selecting an $m$ or a $q$) are different because the two families of models and estimators are different. In particular the families of models have different approximation properties. For instance, the density $2^K\1_{[0,2^{-K})}$ belongs to the two-dimensional model $\overline{S}_{\{1\}}$ but its expansion in the Haar basis has $K$ non-zero coefficients so that it cannot belong to any model $\overline{S}'_q$ with dimension smaller than $K$. As to the estimators, one should notice that histograms are always genuine densities which is not the case of the projection estimators $\widehat{s}_q$.

Penalized projection estimators have been considered by Birg\'e and Massart
(1997) and an improved version is to be found in Chapter~7 of Massart (2007). The method either deals with polynomial collections of linear models, i.e.\ collections for which the number of $D$-dimensional models is bounded by a polynomial in $D$ (which does not apply to our case) or with subset selection within a given basis. Moreover, it requires that $N<n/\log n$ and a bound on $\|\overline{s}_{{\cal I}_{{\cal N}}}\|_\infty$ be known or estimated, as in Section~4.4.4 of Birg\'e and Massart (1997), since the penalty depends on it.

Methods based on wavelet thresholding, as described in Donoho, Johnstone, Kerkyacharian
and Picard (1996) or Kerkyacharian and Picard (2000) (see also the numerous references
therein) typically require the same type of restrictions and, in particular, a known upper bound for $\|s\|_\infty$
in order to properly calibrate the threshold. A noticeable exception appears to be the paper by 
Reynaud-Bouret, Rivoirard and Tuleau-Malot (2011) which is devoted to the estimation of an unknown density on the real line (with possibly unbounded support) by a method which can be viewed as a specific model selection method, the models being linear spaces spanned by finite subsets of a given wavelet basis (for our problem it should be the Haar basis). Their Theorem~1 is some sort of an oracle inequality which does not involve $\|s\|_\infty$ at all but, instead, variance terms which are similar to those in the right-hand side of (\ref{Eq-proj0}). To apply it to densities $s$ belonging to Besov spaces 
$B^\alpha_{p,\infty}(\Bbb{R})$ (with $\alpha>(1/p)-(1/2)$ as is always required) they have to bound these variance terms like we did for (\ref{Eq-proj0}) in Section~\ref{I2b}. They derive risk bounds which show the same dichotomy we mentioned above for histograms. In the ``nice" case (here when $p>2$) the bound does not involve $\|s\|_\infty$. But in the more classical case of 
$p\le2$ they require that $s$ belong to $\Bbb{L}_\infty$ with risk bounds depending on 
$\|s\|_\infty$. 

Aggregation of estimators by selection assumes that preliminary estimators (one for each
model in our case) are given in advance (we should here use the histograms) and typically
leads to a risk bound including a term of the form $n^{-1}\|s\|_\infty\log|{\cal M}|=
n^{-1}N\|s\|_\infty\log2$ so that all such results are useless for $N\ge n$. Moreover,
most of them also require that an upper bound for $\|s\|_\infty$ be known since it enters
the construction of the aggregate estimator. This is the case in Rigollet (2006) (see for
instance his Corollary~2.7) and Juditsky, Rigollet and Tsybakov (2007, Corollary~5.7) since
the parameter $\beta$ that governs their mirror averaging method depends crucially on an
upper bound for $\|s\|_\infty$. As to Samarov and Tsybakov (2005), their Assumption~1
requires that $N$ be not larger than  $C\log n$. Similar restrictions are to be found in Yang
(2000) in his developments for mixing strategies and in Rigollet and Tsybakov (2007) for linear
aggregation of estimators. Lounici (2008) does not assume that $s\in\Bbb{L}_\infty$ but,
instead, that all preliminary estimators are uniformly bounded. One can always truncate the
estimators to get this but, to be efficient, the truncation should be adapted to the unknown
parameter $s$, and therefore chosen from the data in a suitable way. We do not know of any
paper that allows such a data driven choice.

Consequently, none of these results can solve our partition selection problem with arbitrary partitions in a completely satisfactory way when $N$ is at least of size $n$ and whatever the unknown $s\in\Bbb{L}_2(\mu)$. This fact was one motivation for our study of model selection for density estimation with $\Bbb{L}_2$-loss. As already mentioned, some results about adaptive estimation on particular smoothness classes that are akin to model selection with special models and do not assume the boundedness of $s$ can be found in the literature. One can mention estimation of densities belonging to Sobolev classes $W_2^{\alpha}(\Bbb{R})=
B^\alpha_{2,2}(\Bbb{R})$, $\alpha>0$ studied by Efromovich (2008) and densities in 
$B^\alpha_{p,\infty}(\Bbb{R})$, $p>2$ considered by Reynaud-Bouret, Rivoirard and Tuleau-Malot (2011). Both results are quite nice but they address specific situations. The results by Efromovich are actually extremely precise since he does not only get the optimal adaptive rates of convergence but also the exact optimal asymptotic constant which was first computed by Pinsker (1980) for Gaussian ellipsoids. He designs a special estimator of the characteristic function based on an application of the Efromovich-Pinsker method to the empirical characteristic function. Then he proceeds by Fourier inversion. This works remarkably well on Sobolev classes which are defined via the characteristic functions but cannot be extended to more general models. We actually do not know of a {\em general} model selection result that applies to any $s\in\Bbb{L}_2(\mu)$ and arbitrary countable families of  finite-dimensional models (possibly nonlinear). There is a counterpart to this level of generality: our procedure is of a purely abstract nature and not constructive, only indicating what is theoretically feasible. Unfortunately, we are unable to design a practical procedure with similar properties.

\subsection{Some notations\labs{I5}}
Let us now fix our framework and notations. We want to estimate an 
unknown density $s$, with  respect to some {\em probability} measure $\mu$ on the
measurable space $({\cal X}, {\cal W})$, from an i.i.d.\ sample $\bm{X}=(X_1,\ldots,X_n)$ of
random variables $X_i\in{\cal X}$ with distribution $P_s=s\cdot\mu$. The natural domain of application of our results is therefore a compact space ${\cal X}$ with a finite reference measure $\nu$, in which case we shall set $\mu=\nu({\cal X})^{-1}\nu$. Because of this restriction that $\mu$  should be a probability, our result does not apply to estimating densities with respect to the Lebesgue measure on $\Bbb{R}$ but would apply to densities with respect to a Gaussian probability on the line for instance.

Throughout the paper we denote  by $\Bbb{P}_s$ the probability that gives $\bm{X}$ the distribution $P_s^{\otimes n}$ and by $\Bbb{E}_s$ the corresponding expectation operator. For $\Gamma>1$, we  set
\[
\overline{\Bbb{L}}^\Gamma_\infty=
\left\{\left. t\in\overline{\Bbb{L}}_\infty\,\right|\,\|t\|_\infty\le\Gamma\right\}
\]
and, for each $s\in\overline{\Bbb{L}}_2$, we define the function $Q_s$ on $\Bbb{R}_+$ by
\begin{equation}
Q_s(z)=\int[s(x)-z]^2\1_{\{s(x)>z\}}\,d\mu(x)\quad\mbox{for }z\ge0.
\labe{Eq-Qs}
\end{equation}
We measure the performance at $s\in\overline{\Bbb{L}}_2$ of an estimator
$\widehat{s}(\bm{X})\in {\Bbb{L}}_2$ by its quadratic risk
$\Bbb{E}_s\left[d_2^2\left(\widehat{s}(\bm{X}),s\right)\right]$. More generally, if  $(M,d)$ is a
metric space of measurable functions on ${\cal X}$ such that
$M\cap\overline{\Bbb{L}}_1\ne\emptyset$, the quadratic risk of some estimator
$\widehat{s}\in M$ at $s\in M\cap\overline{\Bbb{L}}_1$ is defined as
$\Bbb{E}_s\left[d^2\left(\widehat{s}(\bm{X}),s\right)\right]$. We denote by $|{\cal I}|$ the
cardinality of the set ${\cal I}$ and set $a\vee b$ and $a\wedge b$ for the maximum and the
minimum of $a$ and $b$, respectively. Throughout the paper $C$ (or $C'$, \dots) will denote a universal (numerical) constant and $C(a,b,\ldots)$ or $C_q$ a function of the parameters $a, b,\ldots$ or $q$. Both may vary from line to line. Finally, from now on, {\em countable} will
always mean ``finite or countable".

\Section{Main results\labs{U}}
In order to define estimators based on families of models with bounded metric dimensions,
we shall follow the approach of Birg\'e (2006a) based on what we have called T-estimators. We refer to this paper for the definition and construction of these estimators derived from tests between balls. 

\subsection{Model selection with bounded T-estimators\labs{U1}}
Our first result deals with the performance of special T-estimators that are by construction bounded by $\Gamma$ and therefore belong to $\overline{\Bbb{L}}^\Gamma_\infty$. 
%
\begin{theorem}\lab{T-Main2}
Assume we are given a countable collection $\{\overline{S}_m,m\in{\cal M}\}$ of models in
$\Bbb{L}_2(\mu)$ with metric dimensions bounded respectively by
$\overline{D}_m\ge1/2$ and a family of weights $\Delta_m$ such that
\begin{equation}
\Sigma=1+\sum_{m\in{\cal M}}\exp[-\Delta_m]<+\infty.
\labe{Eq-S8}
\end{equation}
One can build, for each $\Gamma\ge3$, a T-estimator 
$\widehat{s}^\Gamma\in\overline{\Bbb{L}}^\Gamma_\infty$ which  
satisfies, for all $s\in\overline{\Bbb{L}}_2$ and $q\ge1$,
\begin{equation}
\Bbb{E}_s\left[\left\|s-\widehat{s}^\Gamma\right\|^q\right]\le
C_q\Sigma\left[\inf_{m\in{\cal M}}\left\{d_2\!\left(s,\overline{S}_m\right)
+\sqrt{\frac{\Gamma\left(\overline{D}_m\vee\Delta_m\right)}{n}}\right\}
+\sqrt{Q_s(\Gamma)}\right]^q,
\labe{Eq-RB0}
\end{equation}
with $Q_s$ given by (\ref{Eq-Qs}) and $C_q$ some constant depending only on $q$. If
$\|s\|_\infty\le\Gamma$, then 
\begin{equation}
\Bbb{E}_s\left[\left\|s-\widehat{s}^\Gamma\right\|^2\right]\le
C\Sigma\inf_{m\in{\cal M}}\left\{d_2^2\!\left(s,\overline{S}_m\right)
+n^{-1}\Gamma\left(\overline{D}_m\vee\Delta_m\right)\str{4}\!\right\}.
\labe{Eq-RB1}
\end{equation}
\end{theorem}
%
\subsection{General model selection in $\overline{\Bbb{L}}_2$\labs{U2}}
Clearly, the performance of the estimator $\widehat{s}^\Gamma$ provided by Theorem~\ref{T-Main2} depends on the choice of $\Gamma$ since the right-hand side of (\ref{Eq-RB0}) includes a sum of two terms, the first one being increasing with respect to $\Gamma$ and the second one, $[Q_s(\Gamma)]^{1/2}$, nonincreasing. An optimal value of $\Gamma$ should balance between these two terms. Unfortunately both of them depend on the unknown parameter $s$. We therefore need a way to choose $\Gamma$ from the data in order to optimize the bound in (\ref{Eq-RB0}). 

The idea is to build a sequence of estimators $(\widehat{s}^{2^i})_{i\ge2}$ and select a convenient value of $i$ from our data. Since we only have at disposal a single sample 
$\bm{X}$ to build the estimators $\widehat{s}^{2^i}$ and to choose $i$, we shall proceed by sample splitting using one half of the sample for the construction of the estimators and the second half to select a suitable value of $i$. We therefore now consider the general situation where we observe $n=2n'$ i.i.d.\ random variables $X_1,\ldots,X_n$ with an unknown density $s\in\overline{\Bbb{L}}_2$, not necessarily bounded, and have at disposal a countable collection $\{\overline{S}_m,m\in{\cal M}\}$ of models in $\Bbb{L}_2(\mu)$ with metric
dimensions bounded respectively by $\overline{D}_m\ge1/2$ together with a family of weights
$\Delta_m$ which satisfy (\ref{Eq-S8}).  We split our sample $\bm{X}=(X_1,\ldots,X_n)$
into two subsamples $\bm{X}\!_1$ and $\bm{X}\!_2$ of the same size $n'$. We use
$\bm{X}\!_1$ to build the T-estimators $\widehat{s}_i(\bm{X}\!_1)=
\widehat{s}^{2^{i+1}}(\bm{X}\!_1)$, $i\ge1$, which are provided by Theorem~\ref{T-Main2}. 
It then follows from (\ref{Eq-RB0}) that each such estimator satisfies, for $q\ge1$,
\begin{eqnarray*}
\lefteqn{\Bbb{E}_s\left[\left\|s-\widehat{s}_i(\bm{X}\!_1)\right\|^q\right]}\qquad\\&\le&
C_q\Sigma\left\{\inf_{m\in{\cal M}}\left[d_2\!\left(s,\overline{S}_m\right)+
\left(\frac{2^i\left(\overline{D}_m\vee\Delta_m\right)}{n}\right)^{1/2}\right]
+\sqrt{Q_s(2^{i+1})}\right\}^q,
\end{eqnarray*}
with $Q_s$ given by (\ref{Eq-Qs}). We now work conditionally on $\bm{X}\!_1$, fix a
convenient value of $A\ge1$ (for instance $A=1$ if we just want to bound the quadratic
risk) and use the second half of the sample $\bm{X}\!_2$ to select one estimator among
the previous family. This requires a special argument to select a density from an unbounded sequence which is provided by the following proposition to be proved in Section~\ref{T4}.
%
\begin{proposition}\lab{P-Main3}
Let $(t_i)_{i\ge1}$ be a sequence of densities such that 
$t_i\in\overline{\Bbb{L}}_\infty^{2^{i+1}}$ for each $i$ and $\bm{X}$ be an $n$-sample with density $s\in\overline{\Bbb{L}}_2$. Given $A\ge1$, one can design an estimator 
$\widehat{s}_A(\bm{X})$ such that
\[
\Bbb{E}_s\left[d_2^q\!\left(\widehat{s}_A,s\right)\right]\le C(A,q)\inf_{i\ge1}
\left[d_2(s,t_i)\vee\sqrt{n^{-1}i2^i}\right]^q\quad\mbox{for }1\le q<2A/\log2.
\]
\end{proposition}
The selection, based on the sample $\bm{X}\!_2$, of an estimator in the sequence 
$\left(\widehat{s}_i(\bm{X}\!_1)\right)_{i\ge1}$ according to Proposition~\ref{P-Main3} results in a new estimator $\widetilde{s}_A(\bm{X})$ which satisfies
\[
\Bbb{E}_s\left[\left.\st d_2^q\!\left(\widetilde{s}_A(\bm{X}),s\right)\,\right|\,
\bm{X}\!_1\right]\le C(A,q)\inf_{i\ge1}\left[d_2\!\left(s,\widehat{s}_i(\bm{X}\!_1)\right)\vee
\sqrt{n^{-1}i2^i}\right]^q,
\]
provided that $q<2A/\log2$. Integrating with respect to $\bm{X}\!_1$  and using our
previous risk bound gives
\begin{eqnarray*}
\lefteqn{\Bbb{E}_s\left[\left\|s-\widetilde{s}_A(\bm{X})\right\|^q\right]}\\&\le&
\!C(A,q)\inf_{i\ge1}\left\{\Bbb{E}_s\left[\left\|s-\widehat{s}_i(\bm{X}\!_1)\right\|^q
\right]+\left(n^{-1}i2^i\right)^{q/2}\right\}\\&\le&\!C(A,q)\Sigma\inf_{i\ge1}
\left\{\inf_{m\in{\cal M}}\left[d_2^q\left(s,\overline{S}_m\right)+
\left(\frac{2^i\left(\overline{D}_m\vee\Delta_m\vee
i\right)}{n}\right)^{q/2}\right]  +[Q_s(2^{i+1})]^{q/2}\right\}.
\end{eqnarray*}
For $2^i\le z<2^{i+1}$, $\log z\ge i\log2$ and $Q_s(z)\ge Q_s(2^{i+1})$ since $Q_s$ is
nonincreasing. Modifying accordingly the constants in our bounds, we get the main result of
this paper which provides adaptation with respect to both the models and the truncation constant.
%
\begin{theorem}\lab{T-Main0}
Let $\bm{X}=(X_1,\ldots,X_n)$ with $n\ge2$ be an i.i.d.\ sample with unknown
density $s\in\overline{\Bbb{L}}_2$ and $\{\overline{S}_m,m\in{\cal M}\}$ be a
countable collection of models in $\Bbb{L}_2(\mu)$ with metric dimensions bounded
respectively by $\overline{D}_m\ge1/2$. Let $\{\Delta_m,m\in{\cal M}\}$ be a family of
weights which satisfy (\ref{Eq-S8}) and $Q_s(z)$ be given by (\ref{Eq-Qs}). For each 
$A\ge1$, there exists an estimator $\widetilde{s}_A(\bm{X})$ such that, whatever 
$s\in\overline{\Bbb{L}}_2$ and $1\le q<(2A/\log2)$,
\begin{eqnarray}
\lefteqn{\Bbb{E}_s\left[\left\|s-\widetilde{s}_A(\bm{X})\right\|^q\right]}\nonumber\\
&\le&\!\!C(A,q)\Sigma\inf_{z\ge2}\inf_{m\in{\cal M}}\!\left[d_2^q
\left(s,\overline{S}_m\right)+\left(\frac{z\left(\overline{D}_m\vee\Delta_m\vee
\log z\right)}{n}\right)^{q/2}\!+[Q_s(z)]^{q/2}\right]\!.\qquad\;
\labe{Eq-RB4}
\end{eqnarray}
In particular, for $\widetilde{s}=\widetilde{s}_1$ and $s\in\overline{\Bbb{L}}_\infty(\mu)$,
\begin{equation}
\Bbb{E}_s\left[\left\|s-\widetilde{s}(\bm{X})\right\|^2\right]\le C\Sigma
\inf_{m\in{\cal M}}\left[d_2^2\!\left(s,\overline{S}_m\right)+n^{-1}\|s\|_\infty
\left(\overline{D}_m\vee\Delta_m\vee\log\|s\|_\infty\right)\right].
\labe{Eq-RB5}
\end{equation}
\end{theorem}
%

\subsection{Some remarks\labs{U3}}
We see that (\ref{Eq-RB4}) is a generalization of (\ref{Eq-RB0}) and (\ref{Eq-RB5}) of
(\ref{Eq-RB1}) at the modest price of the extra $\log z$ (or $\log\|s\|_\infty$). We do not
know whether this $\log z$ is necessary or not but, in a typical model selection problem, when
$s$ belongs to $\overline{\Bbb{L}}_\infty(\mu)$ but not to $\cup_{m\in{\cal M}}
\overline{S}_m$, the optimal value of $\overline{D}_m$ goes to $+\infty$ with $n$, so that,
for this optimal value, asymptotically 
$\overline{D}_m\vee\Delta_m\vee\log\|s\|_\infty=\overline{D}_m\vee\Delta_m$. 

Up to constants depending on $\|s\|_\infty$, (\ref{Eq-RB5}) is the exact analogue of
(\ref{Eq-M15c}) which shows that, when $s\in\overline{\Bbb{L}}_\infty(\mu)$, all the
results about model selection obtained for the Hellinger distance can be translated in terms
of the $\Bbb{L}_2$-distance.

Note that Theorem~\ref{T-Main0} applies to a single model $\overline{S}$ with metric
dimension bounded by $\overline{D}$, in which case one can use a weight
$\Delta=1/2\le\overline{D}$ which results, if $A=1$, in the risk bound
\begin{equation}
\Bbb{E}_s\left[\left\|s-\widetilde{s}(\bm{X})\right\|^2\right]\le C
\left[d_2^2\left(s,\overline{S}\right)+\inf_{z\ge2}
\left\{\frac{z\left(\overline{D}\vee\log z\right)}{n}+Q_s(z)\right\}\right],
\labe{Eq-RB6}
\end{equation}
and, if $s\in\overline{\Bbb{L}}_\infty(\mu)$,
\begin{equation}
\Bbb{E}_s\left[\left\|s-\widetilde{s}(\bm{X})\right\|^2\right]\le C
\left[d_2^2\!\left(s,\overline{S}\right)+n^{-1}\|s\|_\infty
\left(\overline{D}\vee\log\|s\|_\infty\right)\right].
\labe{Eq-RB7}
\end{equation}
Apart from the extra $\log\|s\|_\infty$, which is harmless when it is smaller than
$\overline{D}$, we recover what we expected, namely the bound (\ref{Eq-M15c}).

Even if $s\in\overline{\Bbb{L}}_\infty(\mu)$ the bound (\ref{Eq-RB4}) may be much
better than (\ref{Eq-RB5}). This is actually already visible with one single model, comparing
(\ref{Eq-RB6}) with (\ref{Eq-RB7}). It is indeed easy to find an example of a very spiky density
$s$ for which (\ref{Eq-RB6}) is much better than (\ref{Eq-RB7}) or the classical bound
(\ref{Eq-proj2}) obtained for projection estimators. Of course, this is just a comparison of
universal bounds, not of the true risk of estimators for a given $s$. 

More surprising is the fact that our estimator can actually dominate a histogram based on the
same model, although our counter-example is rather caricatural and more an advertising
against the use of the $\Bbb{L}_2$-loss than against the use of histogram estimators.
Let us consider a partition ${\cal I}$ of $[0,1]$ into $2D$ intervals $I_j$, $1\le j\le2D$ with
the integer $D$ satisfying $2\le D\le n$ and fix some $\gamma\ge 10$. We then set
$\alpha=\left(\gamma^2n\right)^{-1}$. For $1\le j\le D$, the intervals $I_{2j-1}$ have length
$\alpha$ while the intervals $I_{2j}$ have length $\beta$ with $D(\alpha+\beta)=1$. We
denote by $\overline{S}$ the $2D$-dimensional linear space spanned by the indicator
functions of the $I_j$. It is a model with metric dimension bounded by $D$. We assume that
the underlying density $s$ with respect to Lebesgue measure belongs to $\overline{S}$ and
is defined as
\[
s=p\alpha^{-1}\sum_{j=1}^D\1_{I_{2j-1}}+q\beta^{-1}\sum_{j=1}^D\1_{I_{2j}}
\quad\mbox{with }p=\gamma\alpha\quad\mbox{and}\quad D(p+q)=1,
\]
so that $\beta>q$ since $\alpha<p$. We consider two estimators of $s$ derived from the same
model $\overline{S}$: the histogram $\widehat{s}_{{\cal I}}$ based on the partition ${\cal I}$ and  the
estimator $\widetilde{s}$ based on $\overline{S}$ and provided by Theorem~\ref{T-Main0}.
According to (\ref{Eq-his1}) the risk of $\widehat{s}_{{\cal I}}$ is
\[ 
Dn^{-1}\left[\alpha^{-1}p(1-p)+\beta^{-1}q(1-q)\right]\ge0.9Dn^{-1}\alpha^{-1}p=
0.9D\gamma n^{-1},
\]
since $p\le1/10$. The risk of $\widetilde{s}$ can be bounded by (\ref{Eq-RB4}) with $z=4$
which gives
\[
\Bbb{E}_s\left[\left\|s-\widetilde{s}(\bm{X})\right\|^2\right]\le C\left[4Dn^{-1}+
D\!\int_{I_1}(p/\alpha)^2\,d\mu\right]\\=CD\left[4n^{-1}+p^2\alpha^{-1}\right]=
5CDn^{-1}.
\]
For large enough values of $\gamma$ our estimator is better than the histogram. The
problem actually comes from the observations falling in some of the intervals $I_{2j-1}$ which
will lead to a very bad estimation of $s$ on those intervals. Note that this fact will happen
with a small probability since $Dp=D(\gamma n)^{-1}\le\gamma^{-1}$. Nevertheless,
this event of small probability is important enough to lead to a large risk when we use the
$\Bbb{L}_2$-loss.

\Section{Some applications\labs{X}}

\subsection{Aggregation of preliminary estimators\labs{X4}}
Theorem~\ref{T-Main0} applies in particular to the problem of aggregating preliminary
estimators, built from an independent sample, either by selecting one of them or by combining
them linearily.

\subsubsection{Aggregation by selection\labs{X4a}}
Let us begin with the problem, that we already considered in Section~\ref{U2}, of selecting a
point among a countable family $\{t_m, m\in{\cal M}\}$. Typically, as in Rigollet (2006), the
$t_m$ are preliminary estimators based on an independent sample (derived by sample splitting if necessary) and we want to choose the best one in the family. This is a situation for which one can choose $\overline{D}_m=1/2$ for all $m$ and $A=1$ which leads to the following corollary
%
\begin{corollary}\lab{C-pselect}
Let $\bm{X}=(X_1,\ldots,X_n)$ with $n\ge2$ be an i.i.d.\ sample with unknown
density $s\in\overline{\Bbb{L}}_2$ and $\{t_m,m\in{\cal M}\}$ be a countable
collection of points in $\Bbb{L}_2(\mu)$. Let $\{\Delta_m,m\in{\cal M}\}$ be a
family of weights which satisfy (\ref{Eq-S8}) and $Q_s(z)$ be given by (\ref{Eq-Qs}).
There exists an estimator $\widetilde{s}(\bm{X})$ such that, whatever
$s\in\overline{\Bbb{L}}_2$,
\[
\Bbb{E}_s\left[\left\|s-\widetilde{s}(\bm{X})\right\|^2\right]\le C\Sigma\inf_{z\ge2}
\left\{\inf_{m\in{\cal M}}\left[d_2^2(s,t_m)+
\frac{z(\Delta_m\vee\log z)}{n}\right]+Q_s(z)\right\}.
\]
\end{corollary}
%

\subsubsection{Linear aggregation\labs{X4b}}
Rigollet and Tsybakov (2007) have considered the problem of linear aggregation. Given a finite
set $\{t_1,\ldots, t_N\}$ of preliminary estimators of $s$, they use the observations to build a
linear combination of the $t_j$ in order to get a new and potentially better estimator of $s$. For
$\bm{\lambda}=(\lambda_1,\ldots,\lambda_N)\in\Bbb{R}^N$, let us set
$t_{\bm{\lambda}}=\sum_{j=1}^N\lambda_jt_j$. Rigollet and Tsybakov build a selector
$\widehat{\bm{\lambda}}(X_1,\ldots,X_n)$ such that the corresponding estimator
$\widehat{s}(\bm{X})=t_{\widehat{\bm{\lambda}}}$ satisfies, for all
$s\in\overline{\Bbb{L}}_\infty$,
\begin{equation}
\Bbb{E}_s\left[\left\|s-\widehat{s}(\bm{X})\right\|^2\right]\le
\inf_{\bm{\lambda}\in\Bbb{R}^N}d_2^2\left(s,t_{\bm{\lambda}}\right)+n^{-1}\|s\|_\infty
N.
\labe{Eq-Rig-Tsy}
\end{equation}
Unfortunately, this bound, which is shown to be sharp for such an estimator, can be really poor,
as compared to the minimal risk $\inf_{1\le j\le N}d_2^2(s,t_j)$ of the preliminary estimators
when one of these is already quite good and $n^{-1}\|s\|_\infty N$ is large, which is likely to
happen when $N$ is quite large. Moreover, this result tells nothing when $s\not\in
\overline{\Bbb{L}}_\infty$. In Birg\'e (2006a, Section~9.3) we proposed an alternative way of
selecting a linear combination of the $t_j$ based on T-estimators. In the particular situation of
densities belonging to $\overline{\Bbb{L}}_2$, we proceed as follows: we choose for ${\cal M}$
the collection of all nonvoid subsets $m$ of $\{1,\ldots,N\}$ and, for $m\in{\cal M}$, we take
for $\overline{S}_m$ the linear span of the $t_j$ with $j\in m$ so that the dimension of
$\overline{S}_m$ is bounded by $|m|$ and its metric dimension $\overline{D}_m$ by
$|m|/2$. Since the number of elements of ${\cal M}$ with cardinality $j$ is
$\left(\begin{array}{c}N\\j\end{array}\right)<(eN/j)^j$, we may set $\Delta_m=|m|[2+
\log(N/|m|)]$ so that (\ref{Eq-S8}) is satisfied with $\Sigma<2$. An application of
Theorem~\ref{T-Main0} leads to the following corollary.
%
\begin{corollary}\lab{C-linaggreg}
Let $\bm{X}=(X_1,\ldots,X_n)$ with $n\ge2$ be an i.i.d.\ sample with unknown density
$s\in\overline{\Bbb{L}}_2$ and $\{t_1,\ldots, t_N\}$ be a finite set of points in
$\Bbb{L}_2(\mu)$. Let ${\cal M}$ be the collection of all nonvoid subsets $m$ of
$\{1,\ldots,N\}$ and, for $m\in{\cal M}$, 
\[
\Lambda_m=\left\{\left.\bm{\lambda}\in\Bbb{R}^N\,\right|\,\lambda_j=0
\mbox{ for }j\not\in m\right\}.
\]
For each $A\ge1$, there exists an estimator $\widetilde{s}_A(\bm{X})$ such that, whatever
$s\in\overline{\Bbb{L}}_2$ and $1\le q<(2A/\log2)$,
\[
\Bbb{E}_s\left[\left\|s-\widetilde{s}_A(\bm{X})\right\|^q\right]\le
C(A,q)\inf_{z\ge2}\inf_{m\in{\cal M}}R(q,s,z,m),
\]
where
\[
R(q,s,z,m)=\inf_{\bm{\lambda}\in\Lambda_m}
d_2^q\left(s,t_{\bm{\lambda}}\right)+\left(\frac{z\left[|m|\left(\st1+\log(N/|m|)\right)
\vee\log z\right]}{n}\right)^{q/2}\!+[Q_s(z)]^{q/2}
\]
and $Q_s(z)$ is given by (\ref{Eq-Qs}).
\end{corollary}
There are many differences between this bound and (\ref{Eq-Rig-Tsy}), apart from the nasty
constant $C(A,q)$. Firstly, it applies to densities $s$ that do not belong to
$\overline{\Bbb{L}}_\infty$ and handles the case of $q>2$ for a convenient choice of $A$. Also,
when $s\in\overline{\Bbb{L}}_\infty$ and one of the preliminary estimators is already close to
$s$, it may very well happen, when $N$ is large, that
\[
R\left(2,s,\|s\|_\infty,m\right)\le\inf_{\bm{\lambda}\in\Lambda_m}
d_2^2\left(s,t_{\bm{\lambda}}\right)+n^{-1}\|s\|_\infty
\left[|m|\left(\st1+\log(N/|m|)\right)\vee\log\|s\|_\infty\right]
\]
be much smaller than the right-hand side of (\ref{Eq-Rig-Tsy}) for some $m$ of small
cardinality.

\subsection{Selection of projection estimators\labs{X1}}
In this section, we assume that $s\in\overline{\Bbb{L}}_\infty(\mu)$. This assumption is not
needed for the design of the estimator but only to derive suitable risk bounds. We have at hand
a countable family $\left\{\overline{S}_m, m\in{\cal M}\right\}$ of linear subspaces of
$\Bbb{L}_2(\mu)$ with respective dimensions $D_m$ and we choose corresponding weights
$\Delta_m$ satisfying (\ref{Eq-S8}). For each $m$, we consider the projection estimator
$\widehat{s}_m$ defined in Section~\ref{I2}. Each such estimator has a risk bounded by
(\ref{Eq-proj2}), i.e.
\[
\Bbb{E}_s\left[\|\widehat{s}_m-s\|^2\right]\le\|\overline{s}_m-s\|^2+n^{-1}D_m
\|s\|_\infty,
\]
where $\overline{s}_m$ denotes the orthogonal projection of $s$ onto $\overline{S}_m$. If we
apply Corollary~\ref{C-pselect} to this family of estimators, we get an estimator 
$\widetilde{s}(\bm{X})$ satisfying, for all $s\in\overline{\Bbb{L}}_\infty$,
\[
\Bbb{E}_s\left[\left\|s-\widetilde{s}(\bm{X})\right\|^2\right]\le C\Sigma
\inf_{m\in{\cal M}}\left[\|\overline{s}_m-s\|^2+n^{-1}\|s\|_\infty
\left(D_m\vee\Delta_m\vee\log \|s\|_\infty\right)\right].
\]
With this bound at hand, we can now go back to the problem we considered in 
Section~\ref{R2}, starting with an arbitrary countable family $\{{\cal I}_m, m\in{\cal M}\}$ of finite partitions of ${\cal X}$ and weights $\Delta_m$ satisfying (\ref{Eq-S8}). To each partition ${\cal I}_m$ we associate the linear space $\overline{S}_m$ of piecewise constant functions of the form $\sum_{I\in{\cal I}_m}\beta_I\1_{I}$. The dimension of this linear space is the cardinality of ${\cal I}_m$ and its metric dimension is bounded by $|{\cal I}_m|/2$. If we know that $s\in\overline{\Bbb{L}}_\infty(\mu)$, we can proceed as we just explained, building the
family of histograms $\widehat{s}_{{\cal I}_m}(\bm{X}\!_1)$ corresponding to our partitions and using Corollary~\ref{C-pselect} to get
\begin{equation}
\Bbb{E}_s\left[\left\|s-\widetilde{s}(\bm{X})\right\|^2\right]\le C\Sigma
\inf_{m\in{\cal M}}\left[\|\overline{s}_{{\cal I}_m}-s\|^2+n^{-1}\|s\|_\infty
\left(|{\cal I}_m|\vee\Delta_m\vee\log\|s\|_\infty\right)\right],
\labe{Eq-partsel}
\end{equation}
which should be compared with (\ref{Eq-his5}). Apart from the unavoidable complexity
term $\Delta_m$ due to model selection, we have only lost (up to the universal constant
$C$) the replacement of $|{\cal I}_m|$ by $|{\cal I}_m|\vee\log\|s\|_\infty$. Examples of 
families of partitions and corresponding weights satisfying (\ref{Eq-S8}) are given in Section~9 of Birg\'e (2006a).

In the general case of $s\in\overline{\Bbb{L}}_2(\mu)$, we may apply
Theorem~\ref{T-Main0} to the family of linear models $\left\{\overline{S}_m, m\in{\cal
M}\right\}$ derived from these partitions, getting an estimator $\widetilde{s}$ with a risk satisfying
\[
\Bbb{E}_s\left[\left\|s-\widetilde{s}(\bm{X})\right\|^2\right]\le C\Sigma\inf_{z\ge2}
\left\{\inf_{m\in{\cal M}}\left[\|\overline{s}_{{\cal I}_m}-s\|^2+
\frac{z(|{\cal I}_m|\vee\Delta_m\vee\log z)}{n}\right]+Q_s(z)\right\}.
\]

\subsection{A comparison with Gaussian model selection\labs{X2}}
A benchmark for model selection in general is the particular (simpler) situation of  model
selection for the so-called {\em white noise framework} in which we observe a Gaussian
process $\bm{X}=\{X_z,z\in[0,1]\}$ with $X_z=\int_0^zs(x)\,dx+\sigma W_z$, where $s$ is
an unknown element of  $\Bbb{L}_2( [0,1],dx)$, $\sigma>0$ a known parameter and $W_z$ a
Wiener process. For such a problem, an analogue of Theorem~\ref{T-sevmodel} has been
proved in  Birg\'e (2006a), namely 
%
\begin{theorem}\lab{T-Whitenoise}
Let $\bm{X}$ be the Gaussian process given by 
\[ 
X_z=\int_0^zs(x)\,dx+n^{-1/2}W_z,\quad 0\le z\le1,
\]
where $s$ is an unknown element of  $\Bbb{L}_2( [0,1],dx)$ to be estimated and $W_z$ a
Wiener process. Let $\{\overline{S}_m,m\in{\cal M}\}$ be a countable collection of
models in $\Bbb{L}_2(\mu)$ with metric dimensions bounded respectively by
$\overline{D}_m\ge1/2$. Let $\{\Delta_m,m\in{\cal M}\}$ be a family of weights which
satisfy (\ref{Eq-S8}). There exists an estimator $\widetilde{s}(\bm{X})$ such that, whatever
$s\in\Bbb{L}_2( [0,1],dx)$,
\[
\Bbb{E}_s\left[\left\|s-\widetilde{s}(\bm{X})\right\|^2\right]\le C
\inf_{m\in{\cal M}}\left[d_2^2\!\left(s,\overline{S}_m\right)+n^{-1}
\left(\overline{D}_m\vee\Delta_m\right)\right].
\]
\end{theorem}
Comparing this bound with (\ref{Eq-RB5}) shows that, when
$s\in\overline{\Bbb{L}}_\infty(\mu)$, we get a similar risk bound for estimating the density
$s$ from n i.i.d.\ random variables, apart from an additional factor depending on
$\|s\|_\infty$. Similar analogies are valid with bounds obtained for estimating densities with
squared Hellinger loss or for estimating the intensity of a Poisson process as shown in Birg\'e
(2006a and 2007). Therefore, all the many examples that have been treated in these papers as well as those in Baraud and Birg\'e (2011) could be transferred to the case of density estimation with $\Bbb{L}_2$-loss with minor modifications due to the appearence of 
$\|s\|_\infty$ in the bounds. We leave all these translations as exercices for the concerned reader.

\subsection{Adaptive estimation in Besov spaces\labs{X3}}
The Besov space $B^\alpha_{p,\infty}([0,1])$ with $\alpha,p>0$ is defined in DeVore and
Lorentz (1993) and it is known that a necessary and sufficient condition for
$B^\alpha_{p,\infty}([0,1])\subset\Bbb{L}_2([0,1],dx)$ is $\delta=\alpha+1/2-1/p>0$, which
we shall assume in the sequel. The problem of estimating densities that belong to some Besov 
space $B^\alpha_{p,\infty}([0,1])$ adaptively (i.e.\ without knowing $\alpha$ and $p$) has been 
solved for a long time when $\alpha>1/p$ which is a necessary and sufficient condition for
$B^\alpha_{p,\infty}([0,1])\subset\Bbb{L}_\infty([0,1],dx)$. See for instance Donoho,
Johnstone, Kerkyacharian and Picard (1996), Delyon and Juditsky (1996) (under the assumption that an upper bound for $\|s\|_\infty$ is known) or Birg\'e and Massart
(1997) (with an estimated value of $\|s\|_\infty$). It can be treated in the usual way leading to the minimax rate of convergence $n^{-2\alpha/(2\alpha+1)}$ for the quadratic risk when $n$ goes to infinity. The situation is quite different when $\alpha\le1/p$ even when $\alpha$ and $p$ are known.

\subsubsection{Wavelet expansions\labs{X3y}}
It is known from analysis that functions $s\in\Bbb{L}_2\left([0,1],dx\right)$ can be
represented by their expansion with respect to some orthonormal wavelet basis 
$\{\varphi_{j,k}, j\ge-1, k\in\Lambda(j)\}$ with $|\Lambda(-1)|\le K$ and
$2^j\le|\Lambda(j)|\le K2^j$ for all $j\ge0$. Such a wavelet basis satisfies
\begin{equation}
\left\|\sum_{\:k\in\Lambda(j)\:}|\varphi_{j,k}|\right\|_\infty\le K'2^{j/2}
\;\;\mbox{for }j\ge-1\quad\;\mbox{and}\quad\;
\left\|\sum_{j=-1}^q\sum_{\:k\in\Lambda(j)}\varphi_j^2\right\|_\infty\le K''2^q
\labe{Eq-2z}
\end{equation}
and we can write
\begin{equation}
s=\sum_{j=-1}^\infty\sum_{\:k\in\Lambda(j)\:}\beta_{j,k}\varphi_{j,k},\quad\mbox{with}
\quad\beta_{j,k}=\int\varphi_{j,k}(x)s(x)\,dx.
\labe{Eq-2p}
\end{equation}
Moreover, for a convenient choice of the wavelet basis (depending on $\alpha$), the fact that $s$ belongs to the Besov space $B^\alpha_{p,\infty}([0,1])$ with semi-norm $|s|^\alpha_p$ is equivalent to
\begin{equation}
\sup_{j\ge0}2^{j(\alpha+1/2-1/p)}\left(\sum_{\:k\in\Lambda(j)\:}|
\beta_{j,k}|^p\right)^{1/p}=|s|_{\alpha,p,\infty}<+\infty,
\labe{Eq-7a}
\end{equation}
where $|s|_{\alpha,p,\infty}<+\infty$ is equivalent to the Besov semi-norm $|s|^\alpha_p$.

Moreover, it follows from Birg\'e and Massart (1997 and 2000), as summarized in Birg\'e (2006a,
Proposition~13), that, given the integer $r$, one can find a wavelet basis (depending on $r$) and
a universal family of linear models $\{\overline{S}_m, m\in{\cal M}=\cup_{J\ge0}{\cal M}_J\}$
with respective dimensions $\overline{D}_m$, and weights $\{\Delta_m, m\in{\cal M}\}$
satisfying (\ref{Eq-S8}), with the following properties. Each $\overline{S}_m$ is the linear span
of
$\{\varphi_{-1,k} , k\in\Lambda(-1)\}\cup\{\varphi_{j,k},\, (j,k)\in m\}$ with
$m\subset\cup_{j\ge0}\Lambda(j)$; $\overline{D}_m\vee\Delta_m\le c2^J$ for $m\in{\cal
M}_J$ and
\begin{equation}
\inf_{m\in{\cal M}_J}\inf_{t\in\overline{S}_m}\|s-t\|\le C(\alpha,p)2^{-J\alpha}
|s|_{\alpha,p,\infty}\quad\mbox{for }s\in B^\alpha_{p,\infty}([0,1]),\;\alpha<r.
\labe{Eq-Bes1}
\end{equation}

\subsubsection{The bounded case\labs{X3z}}
Actually, only the assumption that $s\in B^\alpha_{p,\infty}([0,1])\cap
\overline{\Bbb{L}}_\infty(\mu)$, rather than $\alpha>1/p$, is needed to get the optimal rate
of convergence $n^{-2\alpha/(2\alpha+1)}$. Indeed, we may apply the results of
Section~\ref{X1} to the family of models which satisfies (\ref{Eq-Bes1}) and derive an estimator
$\widetilde{s}$ with a risk bounded by
\[
\Bbb{E}_s\left[\left\|s-\widetilde{s}(\bm{X})\right\|^2\right]\le C(\alpha,p)
\inf_{J\ge0}\left[2^{-2J\alpha}\left(|s|_{\alpha,p,\infty}\right)^2+n^{-1}\|s\|_\infty
\left(2^J\vee\log \|s\|_\infty\right)\right].
\]
Choosing $2^J$ of the order of $n^{1/(2\alpha+1)}$ leads to the bound
\[
\Bbb{E}_s\left[\left\|s-\widetilde{s}(\bm{X})\right\|^2\right]\le 
C\left(\alpha,p,|s|_{\alpha,p,\infty},\|s\|_\infty\right)n^{-2\alpha/(2\alpha+1)},
\]
which is valid for all $s\in B^\alpha_{p,\infty}([0,1])\cap\overline{\Bbb{L}}_\infty(\mu)$,
whatever $\alpha<r$ and $p$ and although $\alpha$, $p$, $|s|_{\alpha,p,\infty}$ and
$\|s\|_\infty$ are unknown.

\subsubsection{Further upper bounds for the risk\labs{X3b}}
When $\alpha\le1/p$, i.e.\ $0<\delta\le1/2$, $s$ may be unbounded and the classical 
theory does not apply any more. Results that do not involve $\|s\|_\infty$ are available in Efromovich (2008) for Sobolev classes $W_2^\alpha(\Bbb{R})=B^\alpha_{2,2}(\Bbb{R})$ and for Besov spaces 
$B^\alpha_{p,\infty}(\Bbb{R})$ with $p>2$ in Reynaud-Bouret, Rivoirard and Tuleau-Malot (2011). Nevertheless a general formula for the adaptive minimax risk over balls in 
$B^\alpha_{p,\infty}([0,1])$ for $p\le2$ and $1/p-1/2<\alpha\le1/p$ is presently unknown. Our study will not, unfortunately, solve this problem but, at least, provide some partial information. 
In this section we assume that $\alpha\le1/p$ and restrict ourselves to the case $p\le2$ so that $\delta\le\alpha$. 

We consider the wavelet expansion of $s$ which has been described in Section~\ref{X3y} and, to avoid unnecessary complications, we also assume that $|s|_{\alpha,p,\infty}\ge1$. In what follows, the generic
constant $C$ (changing from line to line) depends on the choice of the basis and $\delta$. Since 
$p\le2$, by (\ref{Eq-7a}),
\[
\left(\sum_{\:k\in\Lambda(j)\:}\beta_{j,k}^2\right)^{1/2}\le
\left(\sum_{\:k\in\Lambda(j)\:}|\beta_{j,k}|^p\right)^{1/p}\le|s|_{\alpha,p,\infty}
2^{-j(\alpha+1/2-1/p)}=|s|_{\alpha,p,\infty}2^{-j\delta},
\]
hence, for $J\in\Bbb{N}$,
\begin{equation}
\left\|\sum_{j>J}\sum_{\:k\in\Lambda(j)\:}\beta_{j,k}\varphi_{j,k}\right\|^2=
\sum_{j>J}\sum_{\:k\in\Lambda(j)\:}\beta_{j,k}^2\le|s|_{\alpha,p,\infty}^2\sum_{j>J}
2^{-2j\delta}=|s|_{\alpha,p,\infty}^22^{-2J\delta}.
\labe{Eq-8b}
\end{equation}
The simplest estimators of $s$ are the projection estimators $\widehat{s}_q$ over the linear spaces
$\overline{S}'_q$ where $\overline{S}'_q$ is spanned by $\{\varphi_{j,k}, -1\le j\le q,
k\in\Lambda(j)\}$
\[
\widehat{s}_q(\bm{X})=\sum_{j=-1}^q\sum_{\:k\in\Lambda(j)\:}\widehat{\beta}_{j,k}(\bm{X})
\varphi_{j,k},\quad\mbox{with}\quad\widehat{\beta}_{j,k}(\bm{X})=n^{-1}\sum_{i=1}^n
\varphi_{j,k}(X_i),
\]
The risk of these estimators can be bounded using (\ref{Eq-proj2}), (\ref{Eq-2z}) and (\ref{Eq-8b}) by
\[
\Bbb{E}_s\left[\left\|s-\widehat{s}_q(\bm{X})\right\|^2\right]\le
d_2^2\left(s,\overline{S}'_q\right)+C2^q/n\le2^{-2q\delta}|s|^2_{\alpha,p,\infty}+C2^q/n.
\]
A convenient choice of $q$, depending on $\delta$ (therefore nonadaptive), then leads to
\[
\Bbb{E}_s\left[\left\|s-\widehat{s}_q(\bm{X})\right\|^2\right]\le
C|s|^2_{\alpha,p,\infty}n^{-2\delta/(2\delta+1)}.
\]
In particular, when $p=2$ we recover the usual minimax rate $n^{-2\alpha/(1+2\alpha)}$ for all values of $\alpha$ but without adaptation. 

One can actually choose $q$ from the data using a penalized least squares estimator and get a similar risk bound without knowing $\delta$ as shown by Theorem~7.5 of Massart (2007) which proves adaptation to the minimax risk when $p=2$. It also leads to an adaptive risk bound for the case 
$\alpha\le1/p$, $p<2$ (hence $\delta<\alpha$), without the restriction $s\in\overline{\Bbb{L}}_\infty([0,1])$ but with a rate which is then slower than $n^{-2\alpha/(1+2\alpha)}$. 

Let us now see what our method can do. Since $s$ is a density, it follows from (\ref{Eq-2p})
and (\ref{Eq-2z}) that $|\beta_{-1,k}|\le\|\varphi_{-1,k}\|_\infty\le K'/\sqrt{2}$, hence
\[
\left\|\sum_{k\in\Lambda(-1)}\beta_{-1,k}\varphi_{-1,k}\right\|_\infty\le
\left(K'/\sqrt{2}\right)\left\|\sum_{k\in\Lambda(-1)}|\varphi_{-1,k}|\right\|_\infty\le
K'^2/2.
\]
Moreover, for $j\ge0$, (\ref{Eq-7a}) implies that $\sup_{\:k\in\Lambda(j)\:}|\beta_{j,k}|
\le2^{-j\delta}|s|_{\alpha,p,\infty}$. Therefore, by (\ref{Eq-2z}),
\[
\left\|\sum_{\:k\in\Lambda(j)\:}\beta_{j,k}\varphi_{j,k}\right\|_\infty\le K'2^{-j(\alpha-1/p)}
|s|_{\alpha,p,\infty},
\]
and, for $J\ge0$,
\[
\left\|\sum_{j=0}^J\sum_{\:k\in\Lambda(j)\:}\beta_{j,k}\varphi_{j,k}\right\|_\infty\le
\left\{\begin{array}{ll}C|s|_{\alpha,p,\infty}&\quad\mbox{if }\alpha>1/p;\\
C(J+1)|s|_{\alpha,p,\infty}&\quad\mbox{if }\alpha=1/p;\\C2^{J(1/p-\alpha)}
|s|_{\alpha,p,\infty}&\quad\mbox{if }\alpha<1/p.\end{array}\right. 
\]
Finally,
\[
\left\|\sum_{j=-1}^J\sum_{\:k\in\Lambda(j)\:}\beta_{j,k}\varphi_{j,k}\right\|_\infty
\!\le C_0L_J|s|_{\alpha,p,\infty}\quad\mbox{ with }\quad L_J=\left\{\begin{array}{ll}1&
\quad\mbox{if }\alpha>1/p;\\(J+1)&\quad\mbox{if }\alpha=1/p;\\
2^{J(1/p-\alpha)}&\quad\mbox{if }\alpha<1/p.\end{array}\right. 
\]
Observing that if $s=u+v$ with $\|u\|_\infty\le z$, then $Q_s(z)\le\|v\|^2$, we can
conclude from (\ref{Eq-8b}) that
\[ 
Q_s\left(C_0L_J|s|_{\alpha,p,\infty}\right)\le2^{-2J\delta}|s|_{\alpha,p,\infty}^2.
\]
Let us now turn back to the family of linear models described in Section~\ref{X3y} that satisfy
(\ref{Eq-Bes1}). Theorem~\ref{T-Main0} asserts the existence of an estimator
$\widetilde{s}(\bm{X})$ based on this family of models and satisfying
\[
\Bbb{E}_s\left[\left\|s-\widetilde{s}(\bm{X})\right\|^2\right]\le C\inf_{z\ge2}
\inf_{m\in{\cal M}}\left[d_2^2\left(s,\overline{S}_m\right)+
\frac{z\left(\overline{D}_m\vee\Delta_m\vee\log z\right)}{n}+Q_s(z)\right].
\]
Given the integers $J,J'$, we may set $z=z_{J'}=C_0L_{J'}|s|_{\alpha,p,\infty}$ and restrict the
minimization to $m\in{\cal M}_J$ which leads by (\ref{Eq-Bes1}) to
\[
\Bbb{E}_s\left[\left\|s-\widetilde{s}(\bm{X})\right\|^2\right]\le C
\left[|s|_{\alpha,p,\infty}^2\left(2^{-2J\alpha}+2^{-2J'\delta}\right)+
n^{-1}L_{J'}|s|_{\alpha,p,\infty}\left(2^J\vee\log z_{J'}\right)\right].
\]
Since $L_{J'}\left(2^J\vee\log z_{J'}\right)$ is a nondecreasing function of both $J$ and $J'$,
this last bound is optimized when $J\alpha$ and $J'\delta$ are approximately equal which
leads to choosing the integer $J'$ so that $J\alpha/\delta\le J'<J\alpha/\delta+1$, hence
$2^{-2J'\delta}\le 2^{-2J\alpha}$. Assuming, moreover, that 
$2^J\ge\log|s|_{\alpha,p,\infty}$, which implies that  $2^J\ge C'\log z_{J'}$, we get
\[
\Bbb{E}_s\left[\left\|s-\widetilde{s}(\bm{X})\right\|^2\right]\le C|s|^2_{\alpha,p,\infty}
\left[2^{-2J\alpha}+2^J\left(n|s|_{\alpha,p,\infty}\right)^{-1}L_{J'}\right].
\]
We finally fix $J$ so that $2^J\ge G>2^{J-1}$, where $G$ is defined below. This choice ensures
that $G\ge\log|s|_{\alpha,p,\infty}$ for $n$ large enough (depending on
$|s|_{\alpha,p,\infty}$), which we assume here.\\
--- If $\alpha>1/p$ we set $G=\left(n|s|_{\alpha,p,\infty}\right)^{1/(2\alpha+1)}$ which
leads to a risk bound of the form
\[
Cn^{-2\alpha/(2\alpha+1)}\left(|s|_{\alpha,p,\infty}\right)^{(2\alpha+2)/(2\alpha+1)}.
\]
--- If $\alpha=1/p$, $L_J'<J\alpha/\delta+2$ and we take
$G=\left(n|s|_{\alpha,p,\infty}/\log n\right)^{1/(2\alpha+1)}$ which leads to the risk bound
\[
C(n/\log n)^{-2\alpha/(2\alpha+1)}
\left(|s|_{\alpha,p,\infty}\right)^{(2\alpha+2)/(2\alpha+1)}.
\]
--- Finally, for $\alpha<1/p$, $L_{J'}<\sqrt{2}\,2^{(J\alpha/\delta)(1/p-\alpha)}$ and we set
$G=\left(n|s|_{\alpha,p,\infty}\right)^{1/[\alpha+1+\alpha/(2\delta)]}$ which
leads to the bound 
\[
Cn^{-2\alpha/[\alpha+1+\alpha/(2\delta)]}
\left(|s|_{\alpha,p,\infty}\right)^{(2+(\alpha/\delta)/[\alpha+1+\alpha/(2\delta)]}.
\]

\subsubsection{Some lower bounds\labs{X3a}}
Lower bounds of the form $n^{-2\alpha/(1+2\alpha)}$ for the minimax risk over Besov balls are
well-known (deriving from lower bounds for H\"older spaces) and they are sharp for
$\alpha>1/p$, as shown in Donoho, Johnstone, Kerkyacharian and Picard (1996). To derive
new lower bounds for the case $\alpha<1/p$ we shall use the following proposition which results easily from classical arguments of Le Cam (1973) --- see also Donoho and Liu (1987) or Yu (1997) ---.
%
\begin{proposition}\lab{P-Mino}
Let $X_1,\ldots,X_n$ be i.i.d.\ observations with an unknown density belonging to a subset 
${\cal S}$ of $\overline{\Bbb{L}}_1(\mu)$ and $d$ a distance on ${\cal S}$. Let $t,u\in{\cal S}$
such that
\[
h(t,u)=h(P_t,P_u)= an^{-1/2},\quad a<2^{-1/2}.
\]
Whatever the estimator $\widehat{s}$ with values in ${\cal S}$ and $p\ge1$,
\begin{equation}
\max\left\{\Bbb{E}_t\left[\st d^p(\widehat{s},t)\right],
\Bbb{E}_u\left[\st d^p(\widehat{s},u)\right]\right\}\ge2^{-p}\left(1-a\sqrt{2}\right)d^p(t,u).
\labe{Eq-A8}
\end{equation}
\end{proposition}
Let us consider some probability density $f\in B^\alpha_{p,\infty}([0,1])$ with compact support included in $(0,1)$ and Besov semi-norm $|f|^\alpha_p$. We set $g(x)=af(2anx)$ for some  $a>(2n)^{-1}$ to be fixed later. Then $g(x)=0$ for $x\not\in\left(0,(2an)^{-1}\right)$,
\[
\|g\|_q=a(2an)^{-1/q}\|f\|_q\qquad\mbox{and}\qquad
|g|^\alpha_p=a(2an)^{\alpha-1/p}|f|^\alpha_p.
\]
Let us now set $t=g+\left[1-(2n)^{-1}\right]\1_{[0,1]}$, so that  $t$ is a density belonging to
$B^\alpha_{p,\infty}([0,1])$ with Besov semi-norm
\[ 
|t|^\alpha_p= |g|^\alpha_p=Ka^{1+\alpha-1/p} n^{\alpha-1/p}\quad\mbox{with }
K=2^{\alpha-1/p}|f|^\alpha_p. 
\] 
For a given value of the constant $K'>0$, the choice
$a=\left[K'n^{1/p-\alpha}\right]^{1/(1+\alpha-1/p)}>(2n)^{-1}$ (at least for $n$ large) leads
to $|t|^\alpha_p= KK'$ so that  $K'$ determines $|t|^\alpha_p$. We also consider the density
$u(x)=t(1-x)$ which has the same Besov semi-norm. Then 
\[ 
h^2(t,u)=\int_0^{(2an)^{-1}}\left(\sqrt{g+\left[1-(2n)^{-1}\right]}-\sqrt{1-(2n)^{-1}}\right)^2
<\int_0^{(2an)^{-1}}\!g=(2n)^{-1}
\] 
and it follows from Proposition~\ref{P-Mino} that any estimator $\widehat{s}$ based on $n$ i.i.d.\ observations
satisfies
\[
\max\left\{\Bbb{E}_t\left[\|t-\widehat{s}\|^2\right],
\Bbb{E}_u\left[\|u-\widehat{s}\|^2\right]\right\}\ge
C\|t-u\|^2=2C\|g\|^2=Can^{-1}\|f\|^2.
\]
Since $an^{-1}=K'^{1/(\delta+1/2)}n^{-2\delta/(\delta+1/2)}$, we finally get
\[
\max\left\{\Bbb{E}_t\left[\|t-\widehat{s}\|^2\right],
\Bbb{E}_u\left[\|u-\widehat{s}\|^2\right]\right\}\ge
C'\left(|t|^\alpha_p\right)^{2/(2\delta+1)}n^{-4\delta/(2\delta+1)},
\]
where $C'$ depends on $K'$, $\|f\|$, $|f|^\alpha_p$ and $\delta$. 

\subsubsection{Conclusion\labs{X3x}}
In the case $\alpha>1/p$, the estimator that we built in Section~\ref{X3b} has the usual rate
of convergence with respect to $n$, namely $n^{-2\alpha/(2\alpha+1)}$, which is known to be
optimal, and we can extend the result to the borderline case $\alpha=1/p$ with only a
logarithmic loss. We do not know whether this additional logarithmic factor is necessary or not. When $\alpha\le1/p$ only partial results are known which do not involve $\|s\|_\infty$. Efromovich (2008) proves the same adaptive estimation rate $n^{-2\alpha/(2\alpha+1)}$ for the Sobolev spaces $W_2^\alpha(\Bbb{R})\subsetneq B^\alpha_{2,\infty}(\Bbb{R})$ (and even gets the exact asymptotic constants) and Reynaud-Bouret, Rivoirard and Tuleau-Malot (2011) for $B^\alpha_{p,\infty}(\Bbb{R})$ with $p>2$. As far as we are aware, nothing is known when $p<2$ and for $B^\alpha_{2,\infty}(\Bbb{R})\setminus W_2^\alpha(\Bbb{R})$. Our lower bound $n^{-4\delta/(2\delta+1)}$ is slower than 
$n^{-2\alpha/(1+2\alpha)}$ when $0<\delta<\alpha[2(\alpha+1)]^{-1}$ or, equivalently, when 
$\alpha+[2(\alpha+1)]^{-1}<1/p$ (which is only possible for $p<2$). This means that the minimax rate $n^{-2\alpha/(1+2\alpha)}$ cannot hold in this range (even without adaptation) but this lower bound tells us nothing when $1/2\le 1/p\le\alpha+[2(\alpha+1)]^{-1}$, in particular when $p=2$.

In the range $p<2$ and $\alpha<1/p$, our upper bound $n^{-2\alpha/[\alpha+1+\alpha/(2\delta)]}$
can be compared with the risk bound for the penalized least squares estimators based on the nested models $\overline{S}'_q$, which is, as we have seen, of order $n^{-2\delta/(2\delta+1)}$. Our rate is better when $\alpha>2\delta/(2\delta+1)$, which is always true for $\alpha\ge1/2$ since
$\delta<1/2$. When $\alpha<1/2$ this requires that $p<2(1-\alpha)/\left(1-2\alpha^2\right)$, which is true independently of $\alpha$ when $p<1+2^{-1/2}$. In any case, these upper bounds never match our lower bound $n^{-4\delta/(2\delta+1)}$ and we have no idea about the true minimax rate (even without adaptation) although we suspect that the rate we have found is suboptimal in the range $\alpha<1/p$.

\subsection{Using a nonlinear model\labs{X7}}
Let us now come back to the parametric problem that we considered in Section~\ref{R9}. We can use the whole set $\overline{S}=\{s_\theta,\, 0<\theta\le1/3\}$ as our model which, in this case, contains the true density $s$ so that there is no approximation term $d_2\left(s_\theta,\overline{S}\right)$. It follows from Proposition~\ref{P-Contrex2} that the dimension of $\overline{S}$ is bounded by 2 so that Theorem~\ref{T-Main0} applies leading to the following bound derived from (\ref{Eq-RB6}):
\begin{equation}
\Bbb{E}_\theta\left[\left\|s_\theta-\widetilde{s}(\bm{X})\right\|^2\right]\le C
\inf_{z\ge2}\left\{n^{-1}z\log z+Q_{s_\theta}(z)\right\}\quad\mbox{for all }\theta\in(0,1/3].
\labe{Eq-Con4}
\end{equation}
For $2\le z<\theta^{-2}$, $Q_{s_\theta}(z)=\theta^3\left(\theta^{-2}-z\right)^2$ and $Q_{s_\theta}(z)=0$ for $z\ge\theta^{-2}$. Optimizing the right-hand side of (\ref{Eq-Con4}) with respect to $z$ leads to the risk bound
\begin{equation}
\Bbb{E}_\theta\left[\left\|s_\theta-\widetilde{s}(\bm{X})\right\|^2\right]\le 
C\theta^{-1}\left[(n\theta)^{-1}\log\left(\theta^{-1}\right)\wedge1\right],
\labe{Eq-Con3}
\end{equation}
which goes to infinity with $\theta^{-1}$.

Let us now see to what extent this result is sharp. It follows from Lemma~\ref{L-Contrex2} that if $\lambda=\theta+(12n)^{-1}$, $h^2(\theta;\lambda)<(8n)^{-1}$, hence $h(\theta,\lambda)<2^{-3/2}n^{-1/2}$. Also
\[
d^2_2(\theta;\lambda)>\theta^{-1}-\left(\theta+(12n)^{-1}\right)^{-1}=
[\theta(n\theta+(1/12))]^{-1}\ge(2\theta)^{-1}\left[(n\theta)^{-1}\wedge12\right].
\]
It then follows from (\ref{Eq-A8}) that, whatever the estimator $\widehat{s}$, we get a lower bound for the risk of the form
\begin{equation}
\max\left\{\Bbb{E}_\theta\left[\left\|s_\theta-\widehat{s}(\bm{X})\right\|^2\right],
\Bbb{E}_\lambda\left[\left\|s_\lambda-\widehat{s}(\bm{X})\right\|^2\right]\right\}\ge
(8\theta)^{-1}\left[(n\theta)^{-1}\wedge12\right],
\labe{Eq-A9}
\end{equation}
which shows that (\ref{Eq-Con3}) is optimal up to the logarithmic factor.

\Section{The construction of T-estimators for $\Bbb{L}_2$-loss\labs{T}}
It will actually require several steps since we cannot simply apply the results of Birg\'e (2006a) straightforwardly. We recall that the construction of T-estimators of parameters belonging to the metric space $(M,d)$ relies on the existence of suitable tests between balls in this space. It is required that the errors of these tests satisfy some specific properties. Unfortunately, in the metric space $\left(\overline{\Bbb{L}}_2,d_2\right)$ tests with such properties cannot exist for arbitrary balls but can be built under the assumption that the centers of the two balls are bounded by some number $\Gamma$, the performance of these tests depending on $\Gamma$. With this result at hand, we can build estimators based on families of special models $S_m$, following Birg\'e (2006a). These models need to be discrete subsets of $\overline{\Bbb{L}}^\Gamma_\infty$ (for some given $\Gamma$) with bounded metric dimension. Since there is no reason that our initial models $\overline{S}_m$ be of this type (think of linear models) we shall have to build such special models $S_m$ satisfying these conditions from ordinary ones. This construction will lead to an estimator $\widehat{s}^\Gamma$ belonging to $\overline{\Bbb{L}}^\Gamma_\infty$, the performance of which is given by Theorem~\ref{T-Main2}. The last step involves the choice of $\Gamma$ among the sequence $(2^{i+1})_{i\ge1}$ as previously explained in Section~\ref{U2}.

\subsection{Tests between $\Bbb{L}_2$-balls \labs{T1}}
To derive such tests, we need a few specific technical tools to deal with the $\Bbb{L}_2$-distance.

\subsubsection{Randomizing our sample\labs{T1a}}
In the sequel we shall make use of randomized tests based on a randomization trick
due to Yang and Barron (1998, page 106) which has the effect of replacing all densities
involved in our problem by new ones which are uniformly bounded away from zero. 
For this, we choose some number $\lambda\in(0,1)$ and consider the mapping $\tau$
from $\overline{\Bbb{L}}_2$ to $\overline{\Bbb{L}}_2$ given by $\tau(u)=\lambda
u+1-\lambda$. Note that $\tau$ is one-to-one and isometric, up to a factor $\lambda$, i.e.\
$d_2(\tau(u),\tau(v))= \lambda d_2(u,v)$. If $u\in\overline{\Bbb{L}}^\Gamma_\infty$,
then $\tau(u)\in\overline{\Bbb{L}}_\infty^{\Gamma'}$ with $\Gamma'=\lambda
\Gamma+1-\lambda$.

Let $s'=\tau(s)$. Given our initial i.i.d.\ sample $\bm{X}$, we want to build new i.i.d.\
variables $X'_1,\ldots,X'_n$ with density $s'$. For this, we consider two independent
$n$-samples, $Z_1,\ldots,Z_n$ and $\varepsilon_1,\ldots,\varepsilon_n$ with respective
distributions $\mu$ and Bernoulli with parameter $\lambda$. Both samples are
independent of $\bm{X}$. We then set $X'_i=\varepsilon_iX_i+(1-\varepsilon_i)Z_i$ for
$1\le i\le n$. It follows that $X'_i$ has density  $s'$ as required. We shall still denote by
$\Bbb{P}_s$ the probability on $\Omega$ that gives $\bm{X}'=(X'_1,\ldots,X'_n)$ the
distribution $P_{s'}^{\otimes n}$. Given two distinct points $t,u\in\overline{\Bbb{L}}_2$
we define a test function $\psi(\bm{X}')$ between $t$ and $u$ as a
measurable function with values in $\{t,u\}$, $\psi(\bm{X}')=t$ meaning deciding $t$ and
$\psi(\bm{X}')=u$ meaning deciding $u$.  

Once we have used the randomization trick of Yang and Barron, for instance with
$\lambda=1/2$, we deal with an i.i.d.\ sample $\bm{X}'$ with a density $s'$ which is bounded
from below by $1/2$ and we may therefore work within the set of densities that satisfy this
property. 

\subsubsection{Preliminary results about tests between some convex sets\labs{T1b}}
The main tool for the design of tests between $\Bbb{L}_2$-balls of densities is the
following proposition which derives from the results of Birg\'e (1984) (keeping here the
notations of that paper) and in particular from Corollary~3.2, specialized to the case of
$I=\{t\}$ and $c=0$.
%
\begin{proposition}\lab{P-test1}
Let ${\cal M}$ be some linear space of finite measures on some measurable space
$(\Omega,{\cal A})$ with a topology of a locally convex separated linear space. Let
${\cal P},{\cal Q}$ be two disjoint sets of probabilities in ${\cal M}$ and $F$ a set of
positive measurable functions on $\Omega$ with the following properties (with respect
to the given topology on ${\cal M}$):

i) ${\cal P}$ and ${\cal Q}$ are convex and compact;

ii) for any $f\in F$ and $0<z<1$ the function $P\mapsto\int f^z\,dP$ is well-defined
and upper semi-continuous on ${\cal P}\cup{\cal Q}$;

iii) for any $P\in{\cal P}$, $Q\in{\cal Q}$, $t\in(0,1)$ and $\varepsilon>0$, there exists
an $f\in F$ such that
\[
(1-t)\int f^t\,dP+t\int f^{1-t}\,dQ<\int (dP)^{1-t}(dQ)^t+\varepsilon;
\]

iv) all probabilities in ${\cal P}$ (respectively in ${\cal Q}$) are mutually absolutely
continuous.\\
Then one can find $\overline{P}\in{\cal P}$ and $\overline{Q}\in{\cal Q}$ such that
\begin{eqnarray*}
\sup_{P\in{\cal P}}\int\left(\frac{\overline{Q}}{\overline{P}}\right)^tdP&=&
\sup_{Q\in{\cal Q}}\int\left(\frac{\overline{P}}{\overline{Q}}\right)^{1-t}dQ\;\;=\;\;
\sup_{P\in{\cal P},Q\in{\cal Q}}\int (dP)^{1-t}(dQ)^t\\&=&
\int\left(d\overline{P}\right)^{1-t}\left(d\overline{Q}\right)^t.
\end{eqnarray*}
\end{proposition}
In Birg\'e (1984) we assumed that ${\cal M}$ was the set of {\em all} finite measures  on
$(\Omega,{\cal A})$ but the proof actually only uses the fact that ${\cal P}$ and ${\cal Q}$
are subsets of ${\cal M}$. Recalling that the Hellinger affinity  between two densities $u$ and $v$
is defined by  $\rho(u,v)=\int\sqrt{uv}\,d\mu=1-h^2(u,v)$, we get the following corollary.
%
\begin{corollary}\lab{C-test2}
Let $\mu$ be a probability measure on $({\cal X},{\cal W})$ and, for $1\le i\le n$, let
$\left({\cal P}_i,{\cal Q}_i\right)$ be a pair of disjoint convex and weakly compact subsets of
$\Bbb{L}_2(\mu)$ such that
\begin{equation}
s>0\;\;\mu\mbox{-a.s.}\qquad\mbox{and}\qquad\int s\,d\mu=1\quad\mbox{for all }
s\in\bigcup_{i=1}^n\left({\cal P}_i\cup{\cal Q}_i\right).
\labe{Eq-test3}
\end{equation}
For each $i$, one can find $p_i\in{\cal P}_i$ and $q_i\in{\cal Q}_i$ such that
\[
\sup_{u\in{\cal P}_i}\int\sqrt{q_i/p_i}\,u\,d\mu=\sup_{v\in{\cal Q}_i}
\int\sqrt{p_i/q_i}\,v\,d\mu=\sup_{u\in{\cal P}_i,v\in{\cal Q}_i}\rho(u,v)
=\rho(p_i,q_i).
\]
Let $\bm{X}=(X_1,\ldots,X_n)$ be a random vector on ${\cal X}^n$ with distribution
$\bigotimes_{i=1}^n(s_i\cdot\mu)$ with $s_i\in{\cal P}_i$ for $1\le i\le n$ and let
$x\in\Bbb{R}$. Then
\[
\Bbb{P}\left[\sum_{i=1}^n\log(q_i/p_i)(X_i)\ge2x\right]\le
e^{-x}\prod_{i=1}^n\rho(p_i,q_i)\le\exp\left[-x-\sum_{i=1}^nh^2(p_i,q_i)\right].
\]
If $\bm{X}$ has distribution $\bigotimes_{i=1}^n(u_i\cdot\mu)$ with $u_i\in{\cal
Q}_i$ for $1\le i\le n$, then
\[
\Bbb{P}\left[\sum_{i=1}^n\log(q_i/p_i)(X_i)\le2x\right]\le
e^x\prod_{i=1}^n\rho(p_i,q_i)\le\exp\left[x-\sum_{i=1}^nh^2(p_i,q_i)\right].
\]
\end{corollary}
%
\noindent{\em Proof:}
We apply Proposition~\ref{P-test1} with $t=1/2$, $({\cal X},{\cal W})=(\Omega,{\cal A})$
and ${\cal M}$ the set of measures of the form $u\cdot\mu$, $u\in\Bbb{L}_2(\mu)$ endowed
with the weak $\Bbb{L}_2$-topology. In view of (\ref{Eq-test3}), ${\cal P}_i$ and ${\cal Q}_i$ can
be identified with two sets of probabilities and we can take for $F$ the set of all positive
functions such that $\log f$ is bounded. As a consequence, all four assumptions of
Proposition~\ref{P-test1} are satisfied. In order to get iii) we simply take for $f$ a suitably
truncated version of $s/u$ when $P=s\cdot\mu$ and $Q=u\cdot\mu$. As to the
probability bounds they derive from classical exponential inequalities, as for Lemma~7 of
Birg\'e (2006a).\cqfd
%

\subsubsection{ Abstract tests between $\Bbb{L}_2$-balls\labs{T1c}}
The purpose of  this section is to prove the following result, which is of independent interest, about the performance of some tests between $\Bbb{L}_2$-balls.
%
\begin{theorem}\lab{T-L2tests}
Let $t,u\in\overline{\Bbb{L}}^\Gamma_\infty$ for some $\Gamma\in(1,+\infty)$. For any
$x\in\Bbb{R}$, there exists a test $\psi_{t,u,x}$  between $t$ and $u$, based on the randomized 
sample $\bm{X}'$ defined in Section~\ref{T1a} with $\lambda=\sqrt{64/65}$, which satisfies
\begin{equation}
\sup_{\left\{s\in\overline{\Bbb{L}}_2\,|\,d_2(s,t)\le
d_2(t,u)/4\right\}}\;\Bbb{P}_s[\psi_{t,u,x}(\bm{X}')=u]
\le \exp\left[-\frac{n\left(\|t-u\|^2+x\right)}{65\Gamma}\right],
\labe{Eq-L2tests1}
\end{equation}
\begin{equation}
\sup_{\left\{s\in\overline{\Bbb{L}}_2\,|\,d_2(s,u)\le
d_2(t,u)/4\right\}}\;\Bbb{P}_s[\psi_{t,u,x}(\bm{X}')=t]
\le \exp\left[-\frac{n\left(\|t-u\|^2-x\right)}{65\Gamma}\right].
\labe{Eq-L2tests2}
\end{equation}
\end{theorem}
%
\noindent{\em Proof:}
It requires several steps. To begin with, we use the randomization trick of  Yang and Barron
described in Section~\ref{T1a}, replacing our original sample $\bm{X}$ by the randomized
sample $\bm{X}'=(X'_1,\ldots,X'_n)$ for some convenient value of $\lambda$ to be chosen
later. Each $X'_i$ has density  $s'\ge1-\lambda$ when $X_i$ has density $s$. Then we build
a test between $t'=\tau(t)$ and $u'=\tau(u)$ based on $\bm{X}'$ and
Corollary~\ref{C-test2}. To do this, we set $\Delta=\|t-u\|$,
\[ 
{\cal P}=\tau\left(\st{\cal B}_{d_2}(t,\Delta/4)\cap\overline{\Bbb{L}}_2\right)
\qquad\mbox{and}\qquad
{\cal Q}=\tau\left(\st{\cal B}_{d_2}(u,\Delta/4)\cap\overline{\Bbb{L}}_2\right).
\]
Then ${\cal P}$ is the subset of the ball ${\cal B}_{d_2}(t',\lambda\Delta/4)$ of those
densities bounded from below by $1-\lambda$. Densities $v$ with such properties are characterized by the fact that $\scal{v}{\1_{\cal X}}=1$ ($\1_{\cal X}\in\Bbb{L}_2(\mu)$ because $\mu$ is a probability) and $\scal{v}{\1_A}\ge(1-\lambda)\mu(A)$ for any measurable set $A$, a fact which is preserved under weak convergence and convex combinations. This shows that ${\cal P}$ is convex and weakly closed. Since 
${\cal B}_{d_2}(t',\lambda\Delta/4)$ is weakly  compact, it is also the case for
${\cal P}$ and the same argument shows that ${\cal Q}$ is also convex and weakly compact.
Moreover $d_2({\cal P},{\cal Q})\ge\lambda\Delta/2$. It then follows from 
Corollary~\ref{C-test2} that one can find $\bar{t}\in{\cal P}$ and $\bar{u}\in{\cal Q}$ such that
\begin{equation}
\Bbb{P}_s\left[\sum_{i=1}^n\log\left(\bar{u}(X'_i)/\bar{t}(X'_i)\right)\ge2y\right]
\le\exp\left[-nh¨^2\left(\bar{t},\bar{u}\right)-y\right]\quad\mbox{if }s\in{\cal P},
\labe{Eq-L2t1}
\end{equation}
while
\begin{equation}
\Bbb{P}_s\left[\sum_{i=1}^n\log\left(\bar{u}(X'_i)/\bar{t}(X'_i)\right)\le2y\right]
\le\exp\left[-nh¨^2\left(\bar{t},\bar{u}\right)+y\right]\quad\mbox{if }s\in{\cal Q}.
\labe{Eq-L2t2}
\end{equation}
Fixing $y=nx/(65\Gamma)$, we finally define $\psi_{t,u,x}(\bm{X}')$ by setting
$\psi_{t,u,x}(\bm{X}')=u$ if and only if
$\sum_{i=1}^n\log\left(\bar{u}(X'_i)/\bar{t}(X'_i)\right)\ge2y$. Since $s'\in{\cal P}$ is
equivalent to $s\in{\cal B}_{d_2}(t,\Delta/4)$ or $d_2(s,t)\le\Delta/4$ and similarily 
$s\in{\cal Q}$ is equivalent to $d_2(s,u)\le\Delta/4$, to derive (\ref{Eq-L2tests1}) and
(\ref{Eq-L2tests2}) from (\ref{Eq-L2t1}) and (\ref{Eq-L2t2}), we just have to show that
$h¨^2\left(\bar{t},\bar{u}\right)\ge(65\Gamma)^{-1}\Delta^2$. We start from the fact, to be
proved below, that 
\begin{equation}
\|\bar{t}\vee\bar{u}\|_\infty\le2(\lambda\Gamma+1-\lambda).
\labe{Eq-auxB1}
\end{equation}
It implies that
\begin{eqnarray*}
h^2\left(\bar{t},\bar{u}\right)&=&\frac{1}{2}\int\left(\sqrt{\bar{t}}-\sqrt{\bar{u}}\right)^2
\,d\mu\;\;=\;\;\frac{1}{2}\int\frac{\left(\bar{t}-\bar{u}\right)^2}{\left(\sqrt{\bar{t}}+
\sqrt{\bar{u}}\right)^2}\,d\mu\;\;\ge\;\;\frac{\|\bar{t}-\bar{u}\|^2}
{16(\lambda\Gamma+1-\lambda)}\\&\ge&\frac{(\lambda\Delta)^2}
{64(\lambda\Gamma+1-\lambda)}\;\;=\;\;\frac{\Delta^2}
{65\Gamma[\lambda+\Gamma^{-1}(1-\lambda)]}\;\;\ge\;\;\frac{\Delta^2}{65\Gamma},
\end{eqnarray*}
since $\Gamma>1$. As to (\ref{Eq-auxB1}), it is a consequence of the next lemma to be proved 
in Section~\ref{P1}. We apply this lemma to the pair $t',u'$ which satisfies 
$\|t'\vee u'\|_\infty\le\lambda\Gamma+1-\lambda$. If (\ref{Eq-auxB1}) were wrong, we could find 
$\bar{t}'\in{\cal P}$ and $\bar{u}'\in{\cal Q}$ with $h\left(\bar{t}',\bar{u}'\right)<h\left(\bar{t},\bar{u}\right)$, which, by Corollary~\ref{C-test2}, is impossible.\cqfd
%
\begin{lemma}\lab{L-l2b}
Let us consider four elements $t,u,v_1,v_2$ in $\overline{\Bbb{L}}_2$ with $t\ne u$,
$v_1\ne v_2$ and $\|t\vee u\|_\infty=B$. If $\|v_1\vee v_2\|_\infty>2B$, there  exists
$v'_1,v'_2\in\overline{\Bbb{L}}_2$ with $d_2(v'_1,t)\le d_2(v_1,t)$, $d_2(v'_2,u)\le
d_2(v_2,u)$ and $h(v'_1,v'_2)<h(v_1,v_2)$.
\end{lemma}
%

\subsection{The performance of T-estimators for discrete models\labs{T2}}
We are now in a position to prove an analogue of Corollary~6 of Birg\'e (2006a).
%
\begin{theorem}\lab{T-Main1}
Assume that we observe $n$ i.i.d.\ random variables with unknown density
$s\in\left(\overline{\Bbb{L}}_2,d_2\right)$ and that we have at disposal a countable family of
discrete subsets $\{S_m\}_{m\in{\cal M}}$ of $\overline{\Bbb{L}}^\Gamma_\infty$ for
some given $\Gamma>1$. Let each set $S_m$ satisfy
\begin{equation}
|S_m\cap{\cal B}_{d_2}(t,x\eta_m)|\le\exp\left[D_mx^2\right]\quad\mbox{for
all }x\ge2\mbox{ and }t\in\overline{\Bbb{L}}_2,
\labe{Eq-M19}
\end{equation}
with $\eta_m>0$, $D_m\ge1/2$,
\begin{equation}
\eta_m^2\ge\frac{273\Gamma D_m}{n}\quad\mbox{for all }m\in{\cal M},\quad\;
\sum_{m\in{\cal M}}\exp\left[-\frac{n\eta_m^2}{1365\Gamma}\right]=\Sigma'<+\infty.
\labe{Eq-S9}
\end{equation}
Then one can build a T-estimator $\widehat{s}$ such that, for all $s\in\overline{\Bbb{L}}_2$,
\begin{equation}
\Bbb{E}_s\left[d_2^q(s,\widehat{s})\st\right]\le C_q(\Sigma'+1)\inf_{m\in{\cal M}}
\left\{d_2(s,S_m)\vee\eta_m\st\right\}^q,\quad\mbox{for all }q\ge1.
\labe{Eq-ribou1}
\end{equation}
\end{theorem}
%
\noindent{\em Proof:}
Since (\ref{Eq-ribou1}) is merely a version of (7.6) of Birg\'e (2006a) with $d=d_2$, we just have
to show that Theorem~5 of this paper applies to our situation. It relies on Assumptions~1 and 3
of the paper. Assumption~3 follows from (\ref{Eq-M19}). As to Assumption~1 (with
$a=n/(65\Gamma)$, $B=B'=1$ and $\delta=4d_2$, hence $\kappa=4$), it is a consequence of
our Theorem~\ref{T-L2tests}. The conditions (7.2) and (7.4) of Birg\'e (2006a) on $\eta_m$ and
$D_m$ follow from (\ref{Eq-S9}).\cqfd\vspace{2mm}\\
In the case of a single $D$-dimensional model $\overline{S}\subset
\overline{\Bbb{L}}^\Gamma_\infty$ we get the following corollary:
%
\begin{corollary}\lab{C-Main0}
Assume that we observe $n$ i.i.d.\ random variables with unknown distribution $P_s$,
$s\in\left(\overline{\Bbb{L}}_2,d_2\right)$ and that we have at disposal a $D$-dimensional
model $\overline{S}\subset\overline{\Bbb{L}}^\Gamma_\infty$  for some given
$\Gamma>1$. One can build a T-estimator $\widehat{s}$ such that, for all
$s\in\overline{\Bbb{L}}_2$,
\[
\Bbb{E}_s\left[\|s-\widehat{s}\|^2\right]\le C\left[\inf_{t\in\overline{S}}
d_2^2(s,t)+n^{-1}D\Gamma\right].
\]
\end{corollary}
%
\noindent{\em Proof:}
By Definition~\ref{D-Bmdim} and the remark following it, for each $\eta_0>0$, one can find an
$\eta_0$-net $S_0\subset\overline{S}$ for $\overline{S}$, hence $S_0\subset
\overline{\Bbb{L}}^\Gamma_\infty$, satisfying (\ref{Eq-M19}) with $D_0=25D/4$. Moreover
$d(s,S_0)\le\eta_0+d\left(s,\overline{S}\right)$. Choosing $\eta_0^2=273\times 25\Gamma
D/4$, we may apply Theorem~\ref{T-Main1}. The result then follows from (\ref{Eq-ribou1})
with $q=2$.\cqfd\vspace{2mm}\\
Theorem~\ref{T-Main1} applies in particular to the special situation of each model $S_m$ being
reduced to  a single point $\{t_m\}$ so that we can take $D_m=1/2$ for each $m$. We then get the
following useful corollary.
%
\begin{corollary}\lab{C-Main1}
Assume that we observe $n$ i.i.d.\ random variables with unknown distribution $P_s$,
$s\in\left(\overline{\Bbb{L}}_2,d_2\right)$ and that we have at disposal a countable subset
$S=\{t_m\}_{m\in{\cal M}}$ of $\overline{\Bbb{L}}^\Gamma_\infty$ for some given
$\Gamma>1$. Let $\{\Delta_m\}_{m\in{\cal M}}$ be a family of weights such that
$\Delta_m\ge1/10$ for all $m\in{\cal M}$ satisfying (\ref{Eq-S8}).
We can build a T-estimator $\widehat{s}$ such that, for all $s\in\overline{\Bbb{L}}_2$,
\[
\Bbb{E}_s\left[d_2^q(s,\widehat{s})\st\right]\le C_q\Sigma\inf_{m\in{\cal M}}
\left\{d_2(s,t_m)\vee\sqrt{\Gamma\Delta_m/n}\st\right\}^q\quad\mbox{for all }q\ge1.
\]
\end{corollary}
%
\noindent{\em Proof:} Let us set here $S_m=\{t_m\}$, $D_m=1/2$ and
$\eta_m=37\sqrt{\Gamma\Delta_m/n}$ for $m\in{\cal M}$. One can then check that
(\ref{Eq-M19}) and (\ref{Eq-S9}) are satified so that (\ref{Eq-ribou1}) holds. Our risk bound
follows.\cqfd

\subsection{Model selection with uniformly bounded models\labs{T3}}
At this stage, there is a major difficulty to apply Theorem~\ref{T-Main1} or
Corollary~\ref{C-Main1} which is to  build suitable subsets $S_m$ (or $S$) of
$\overline{\Bbb{L}}^\Gamma_\infty$ from classical approximating sets (models), finite
dimensional linear spaces for instance, that belong to ${\Bbb{L}}_2(\mu)$. We shall now address this problem.

\subsubsection{The projection operator onto $\overline{\Bbb{L}}^\Gamma_\infty$\labs{T3a}}
Our first task is to define a projection operator $\pi_\Gamma$ from $\Bbb{L}_2(\mu)$
onto $\overline{\Bbb{L}}^\Gamma_\infty$ ($\Gamma>1$) and to study its properties. In
the sequel, we systematically identify a real number $a$ with the function $a\1_{{\cal X}}$
for the sake of simplicity. The following proposition is the corrected version, by Yannick
Baraud, of the initially mistaken result of the author.
%
\begin{proposition}\lab{P-projop}
For $t\in\Bbb{L}_2(\mu)$ and $1<\Gamma<+\infty$ we set 
$\pi_\Gamma(t)=[(t+\gamma)\vee0]\wedge\Gamma$ where $\gamma$ is defined by
$\int[(t+\gamma)\vee0]\wedge\Gamma\,d\mu=1$. Then $\pi_\Gamma$ is the projection
operator from $\Bbb{L}_2(\mu)$ onto the convex set $\overline{\Bbb{L}}^\Gamma_\infty$.
Moreover, if $s\in\overline{\Bbb{L}}_2$ and $\Gamma>2$, then
\[
\|s-\pi_\Gamma(s)\|^2\le\frac{\Gamma^2-\Gamma-1}{\Gamma(\Gamma-2)}Q_s(\Gamma),
\]
with $Q_s(z)$ given by (\ref{Eq-Qs}).
\end{proposition}
%
\noindent{\em Proof:}
First note that the existence of $\gamma$ follows from the continuity and monotonicity of the
mapping $z\mapsto\int[(t+z)\vee0]\wedge\Gamma\,d\mu$ and that $\pi_\Gamma(t)
\in\overline{\Bbb{L}}^\Gamma_\infty$. Since $\overline{\Bbb{L}}^\Gamma_\infty$ is a 
closed convex subset of a Hilbert space, the projection operator $\pi$ onto
$\overline{\Bbb{L}}^\Gamma_\infty$ exists and is characterized by the fact that
\begin{equation}
\scal{t-\pi(t)}{u-\pi(t)}\le0\quad\mbox{for all }u\in\overline{\Bbb{L}}^\Gamma_\infty.
\labe{Eq-pro1}
\end{equation}
Since $\int[u-\pi(t)]\,d\mu=0$ for $u\in\overline{\Bbb{L}}^\Gamma_\infty$,
(\ref{Eq-pro1}) implies that $\scal{t+z-\pi(t)}{u-\pi(t)}\le0$ for $z\in\Bbb{R}$, hence
$\pi(t)=\pi(t+z)$. Since $\pi_\Gamma(t)=\pi_\Gamma(t+z)$ as well, we may assume that
$\int[t\vee0]\wedge\Gamma\,d\mu=1$, hence $\pi_\Gamma(t)=
[t\vee0]\wedge\Gamma$ and $\pi_\Gamma(t)=t$ on the set $0\le t\le\Gamma$. Then, for
$u\in\overline{\Bbb{L}}^\Gamma_\infty$,
\[
\scal{t-\pi_\Gamma(t)}{u-\pi_\Gamma(t)}=\int_{t<0}tu\,d\mu+\int_{t>\Gamma}
(t-\Gamma)(u-\Gamma)\,d\mu\le0,
\]
since $0\le u\le\Gamma$. This concludes the proof that $\pi=\pi_\Gamma$.

Let us now bound $\left\|s-\pi_\Gamma(s)\right\|$ when $s\in\overline{\Bbb{L}}_2$,
setting $s=s\wedge\Gamma+v$ with $v=(s-\Gamma)\1_{s>\Gamma}$. Since there is
nothing to prove when $\|s\|_\infty\le\Gamma$, we assume that $\int v\,d\mu>0$. By
Cauchy-Schwarz Inequality,
\begin{equation}
\left(\int v\,d\mu\right)^2\le\mu(\{s>\Gamma\})\int
v^2\,d\mu\le\Gamma^{-1}\|v\|^2.
\labe{Eq-pro2}
\end{equation}
Moreover, since $\int s\wedge\Gamma\,d\mu<1$, $\pi_\Gamma(s)=
(s+\gamma)\wedge\Gamma$ with $0<\gamma\le1$. Hence
\begin{eqnarray*}
1&=&\int[(s+\gamma)\wedge\Gamma]\,d\mu\;\;\ge\;\;\int(s\wedge\Gamma)\,d\mu
+\gamma\mu(\{s\le\Gamma-\gamma\})\\&\ge&1-\int v\,d\mu+\gamma\left(1-
\frac{1}{\Gamma-\gamma}\right)\;\;>\;\;1-\int v\,d\mu+\gamma
\frac{\Gamma-2}{\Gamma-1}
\end{eqnarray*}
and $\gamma<(\Gamma-1)/(\Gamma-2)\int v\,d\mu$. Now, since $0\le\pi_\Gamma(s)-s
\le\gamma$ when $s\le\Gamma$,
\begin{eqnarray*}
\|s-\pi_\Gamma(s)\|^2&=&\int_{s\le\Gamma}[\pi_\Gamma(s)-s]^2\,d\mu+
\|v\|^2\;\;\le\;\;\gamma\int_{s\le\Gamma}[\pi_\Gamma(s)-s]\,d\mu+\|v\|^2\\
&<&\frac{\Gamma-1}{\Gamma-2}\left(\int v\,d\mu\right)
\int_{s>\Gamma}[s-\pi_\Gamma(s)]\,d\mu+\|v\|^2\\&\le&
\frac{\Gamma-1}{\Gamma-2}\left(\int v\,d\mu\right)^2+\|v\|^2\;\;\le\;\;
\left(1+\frac{\Gamma-1}{\Gamma(\Gamma-2)}\right)\|v\|^2,
\end{eqnarray*}
where we used (\ref{Eq-pro2}). This concludes our proof.\cqfd
%
 
\subsubsection{Selection with uniformly bounded models\labs{T3c}}
Typical models $\overline{S}$ for density estimation in $\Bbb{L}_2(\mu)$ are
finite-dimensional linear spaces which are not subsets of $\overline{\Bbb{L}}^\Gamma_\infty$
but merely spaces of functions with nice approximation properties. To apply
Theorem~\ref{T-Main1} we have to replace them by discrete subsets of
$\overline{\Bbb{L}}^\Gamma_\infty$ that satisfy (\ref{Eq-M19}). Unfortunately, they cannot
simply be derived by a discretization of $\overline{S}$ followed by a projection $\pi_\Gamma$
or a discretization of $\pi_\Gamma\left(\overline{S}\right)$. A more complicated construction
is required to preserve both the metric and approximation properties of $\overline{S}$. It is
provided by the following preliminary result.
%
\begin{proposition}\lab{P-newmodels}
Let $\overline{S}$ be a subset of $\Bbb{L}_2(\mu)$ with metric dimension bounded by $D$.
For $\Gamma>2$ and $\eta>0$, one can find a discrete subset $S'$ of
$\overline{\Bbb{L}}^\Gamma_\infty$ with the following properties:
\begin{equation}
|S'\cap{\cal B}_{d_2}(t,x\eta)|\le\exp\left[9Dx^2\right]\quad\mbox{for
all }x\ge2\mbox{ and }t\in\Bbb{L}_2(\mu);
\labe{Eq-proj11}
\end{equation}
for any $s\in\overline{\Bbb{L}}_2$, one can find some $s'\in S'$ such that
\begin{equation}
\|s-s'\|\le3.1\left[\eta+\inf_{t\in\overline{S}}\|s-t\|\right]+
4.1\left(\frac{\Gamma^2-\Gamma-1}{\Gamma(\Gamma-2)}Q_s(\Gamma)\right)^{1/2}.
\labe{Eq-proj6}
\end{equation}
\end{proposition}
%
\noindent{\em Proof:}
According to Definition~\ref{D-Bmdim}, we choose some $\eta$-net $S_\eta$ for
$\overline{S}$ such that (\ref{Eq-M49}) holds for all $t\in\Bbb{L}_2(\mu)$. Since, by
Proposition~\ref{P-projop}, the operator $\pi_\Gamma$ from $\Bbb{L}_2(\mu)$ to
$\overline{\Bbb{L}}^\Gamma_\infty$ satisfies $\|u-\pi_\Gamma(t)\|\le\|u-t\|$ for
all $u\in\overline{\Bbb{L}}^\Gamma_\infty$, we may apply Proposition~12 of Birg\'e
(2006a) with $M'=\Bbb{L}_2(\mu)$, $d=d_2$, $M_0=\overline{\Bbb{L}}^\Gamma_\infty$,
$T=S_\eta$, $\overline{\pi}=\pi_\Gamma$ and $\lambda=1$. It shows that one can
find a subset $S'$ of $\pi_\Gamma(S_\eta)$ such that (\ref{Eq-proj11}) holds and
$d_2(u,S')\le3.1d_2(u,S_\eta)$ for all $u\in\overline{\Bbb{L}}^\Gamma_\infty$. If $s$
is an arbitrary element of $\overline{\Bbb{L}}_2$, then
\[
d_2\left(\pi_\Gamma(s), S'\right)\le3.1d_2\left(\pi_\Gamma(s), S_\eta\right)\le
3.1\left[d_2\left(\pi_\Gamma(s), s\right)+d_2\left(s,\overline{S}\right)+\eta\right],
\]
hence
\begin{equation}
d_2\left(s, S'\right)\le3.1\left[d_2\left(s,\overline{S}\right)+\eta\right]+
4.1d_2\left(\pi_\Gamma(s), s\right).
\labe{Eq-th5a}
\end{equation}
The conclusion follows from Proposition~\ref{P-projop}.\cqfd\vspace{2mm}\\
We are now in a position to derive our main result about bounded model selection. We
start with a countable collection $\{\overline{S}_m,m\in{\cal M}\}$ of models in
$\Bbb{L}_2(\mu)$ with metric dimensions bounded respectively by
$\overline{D}_m\ge1/2$ and a family of weights $\Delta_m$ satisfying (\ref{Eq-S8}).
We  fix some $\Gamma\ge3$ and, for each $m\in{\cal M}$, we set
\[
\eta_m=\left[\left(50\sqrt{\overline{D}_m}\right)
\vee\left(37\sqrt{\Delta_m}\right)\right]\sqrt{\Gamma/n}.
\]
By Proposition~\ref{P-newmodels} (with $\eta=\eta_m$), each $\overline{S}_m$ gives
rise to a subset $S^\Gamma_m$ which satisfies (\ref{Eq-M19}) with $D_m=9\overline{D}_m$.
It follows from our choice of $\eta_m$ that (\ref{Eq-S9}) is also satisfied so that we may apply
Theorem~\ref{T-Main1} to the family of sets $\left\{S^\Gamma_m,m\in{\cal M}\right\}$.
This results in a T-estimator $\widehat{s}^\Gamma$ such that, for all $s\in \overline{\Bbb{L}}_2$,
\[
\Bbb{E}_s\left[d_2^q\!\left(s,\widehat{s}^\Gamma\right)\right]\le C_q\Sigma
\inf_{m\in{\cal M}}\left\{d_2\!\left(s,S^\Gamma_m\right)\vee\eta_m
\st\right\}^q\quad\mbox{for }q\ge1.
\]
We also derive from Proposition~\ref{P-newmodels} that
\[
d_2\!\left(s,S^\Gamma_m\right)\le3.1\left[\eta_m+\inf_{t\in\overline{S}_m}
\|s-t\|\right]+4.1\sqrt{(5/3)Q_s(\Gamma)}.
\]
Putting the bounds together and rearranging the terms leads to Theorem~\ref{T-Main2}.
%

\subsection{An additional selection theorem\labs{T4}}
In order to derive Theorem~\ref{T-Main0} we need an additional selection step in order to
choose a proper estimator among the sequence of estimators $(\widehat{s}^{2^i})_{i\ge1}$.
We start with a general selection result, to be proved in Section~\ref{P3}, that we state
for an arbitrary statistical framework since it may apply to other situations than density
estimation from an i.i.d.\ sample. We observe some random object $\bm{X}$ with distribution
$P_s$ on ${\cal X}$ where $s$ belongs to a metric space $M$ (carrying a distance $d$) which indexes a family ${\cal P}=\{P_t, t\in M\}$ of probabilities on ${\cal X}$. 
%
\begin{theorem}\lab{T-Selection}
Let $(t_p)_{p\ge1}$ be a sequence in $M$ such that the following assumption holds: for all
pairs $(n,p)$ with $1\le n<p$ and all $x\in\Bbb{R}$, one can find a test $\psi_{t_n,t_p,x}$
based on the observation $\bm{X}$ and satisfying
\begin{equation}
\sup_{\{s\in M\,|\,d(s,t_n)\le d(t_n,t_p)/4\}}\;\Bbb{P}_s[\psi_{t_n,t_p,x}(\bm{X})=t_p]
\le B\exp\left[-a2^{-p}d^2(t_n,t_p)-x\right];
\labe{Eq-tes1}
\end{equation}
\begin{equation}
\sup_{\{s\in M\,|\,d(s,t_p)\le d(t_n,t_p)/4\}}\;\Bbb{P}_s[\psi_{t_n,t_p,x}(\bm{X})=t_n]
\le B\exp\left[-a2^{-p}d^2(t_n,t_p)+x\right];
\labe{Eq-tes2}
\end{equation}
with positive constants $a$ and $B$ independent of $n,p$ and $x$. For each $A\ge1$, one can
design an estimator $\widehat{s}_A$ with values in $\{t_p,\,p\ge1\}$ such that, for all $s\in M$,
\begin{equation}
\Bbb{E}_s\left[d^q\left(\widehat{s}_A,s\right)\right]\le BC(A,q)\inf_{p\ge1}
\left[d(s,t_p)\vee\sqrt{a^{-1}p2^p}\right]^q\quad\mbox{for }1\le q<2A/\log2.
\labe{Eq-RB3}
\end{equation}
\end{theorem}
This general result applies to our specific framework of density estimation based on an
observation $\bm{X}$ with distribution $P_s$, $s\in\overline{\Bbb{L}}_2$, provided
that the sequence $(t_p)_{p\ge1}$ be suitably chosen. We shall simply assume here that
$t_p\in\overline{\Bbb{L}}_2$ with $\|t_p\|_\infty\le 2^{p+1}$ for each $p\ge1$. This implies
that, for $1\le i<j$, $t_i$ and $t_j$ belong to $\overline{\Bbb{L}}^{2^{j+1}}_\infty$ so that
Theorem~\ref{T-L2tests} applies with $\bm{X}$ replaced by the randomized sample
$\bm{X}'$ and the assumption of Theorem~\ref{T-Selection} is therefore satisfied with
$d=d_2$, $B=1$ and $a=n/65$, leading to Proposition~\ref{P-Main3}.
%

\Section{Proofs\labs{P}}

\subsection{Proof of Proposition~\ref{P-Contrex2}\labs{P0}}
For simplicity, we shall write $h(\theta,\lambda)$ for $h(s_\theta,s_\lambda)$ and
analogously  $d_2(\theta,\lambda)$ for $d_2(s_\theta,s_\lambda)$ and start with a preliminary lemma.
%
\begin{lemma}\lab{L-Contrex2}
For the parametric problem described in Proposition~\ref{P-Contrex2}, the following holds for all $\theta$ and $\lambda$ in $(0,1/3]$:
\begin{equation}
h^2(\theta,\lambda)=C(\theta,\lambda)|\theta-\lambda|
\quad\mbox{with }2/9<C(\theta,\lambda)<3/2
\labe{Eq-Con1}
\end{equation}
and
\begin{equation}
d_2^2(\theta,\lambda)=C(\theta,\lambda)\left|\theta^{-1}-\lambda^{-1}\right|
\quad\mbox{with }1<C(\theta,\lambda)<3.
\labe{Eq-Con2}
\end{equation}
\end{lemma}
%
{\em Proof:}
Let us first evaluate $h^2(\theta,\lambda)$ for $0<\theta<\lambda\le1/3$. Setting
$\beta_\theta=\left(\theta^2+\theta+1\right)^{-1}\in[9/13,1)$, we get
\begin{eqnarray*}
2h^2(\theta,\lambda)&=&\int_0^1\left(\sqrt{s_\theta(x)}-\sqrt{s_\lambda(x)}\right)^2dx\\
&=&\theta^{3}\left(\theta^{-1}-\lambda^{-1}\right)^2+\left(\lambda^{3}-\theta^{3}\right)
\left(\lambda^{-1}-\sqrt{\beta_\theta}\right)^2+\left(1-\lambda^{3}\right)
\left(\sqrt{\beta_\theta}-\sqrt{\beta_\lambda}\right)^2\\&=&(\lambda-\theta)
\frac{\theta}{\lambda}\left(1-\frac{\theta}{\lambda}\right)
+(\lambda-\theta)\left[1+\frac{\theta}{\lambda}+\left(\frac{\theta}{\lambda}\right)^2
\right]\left(1-\lambda\sqrt{\beta_\theta}\right)^2\\&&\mbox{}+
\left(1-\lambda^{3}\right)\left(\sqrt{\beta_\theta}-\sqrt{\beta_\lambda}\right)^2.
\end{eqnarray*}
Note that the monotonicity of $\theta\mapsto\beta_\theta$ implies that
\[
4/9<\left(1-\lambda\sqrt{\beta_\theta}\right)^2<1,\qquad
\sqrt{\beta_\theta}+\sqrt{\beta_\lambda}>2\sqrt{\beta_{1/3}}=6/\sqrt{13}
\]
and
\begin{equation}
0<\beta_\theta-\beta_\lambda=\frac{(\lambda-\theta)(\lambda+\theta+1)}
{(\theta^2+\theta+1)(\lambda^2+\lambda+1)}<\lambda-\theta.
\labe{Es1}
\end{equation}
It follows that 
\[
0<\left(\sqrt{\beta_\theta}-\sqrt{\beta_\lambda}\right)^2=
\frac{(\beta_\theta-\beta_\lambda)^2}{\left(\sqrt{\beta_\theta}+\sqrt{\beta_\lambda}
\right)^2}<\frac{13}{36}(\lambda-\theta)^2=\frac{13\lambda}{36}
(\lambda-\theta)\left(1-\frac{\theta}{\lambda}\right)
\]
and
\[
0<\left(1-\lambda^{3}\right)\!\left(\sqrt{\beta_\theta}-\sqrt{\beta_\lambda}\right)^2\!
<\frac{13\lambda\left(1-\lambda^{3}\right)}{36}(\lambda-\theta)\!
\left(1-\frac{\theta}{\lambda}\right)<
\frac{2(\lambda-\theta)}{17}\!\left(1-\frac{\theta}{\lambda}\right).
\]
We can therefore write 
\[ 
G=2(\lambda-\theta)^{-1}h^2(\theta,\lambda)=z(1-z)+
c_1(\theta,\lambda)\left(1+z+z^2\right)+c_2(\theta,\lambda)(1-z),
\]
with $z=\theta/\lambda\in(0,1)$, $4/9<c_1(\theta,\lambda)<1$ and
$0<c_2(\theta,\lambda)<2/17$. Since, for given values of $c_1$ and $c_2$, the right-hand side is increasing with respect to $z$, $4/9<c_1<G<3c_1<3$ and (\ref{Eq-Con1}) follows.

Let us now proceed with the $\Bbb{L}_2$-distance $d_2$.
\begin{eqnarray*}
d_2^2(\theta,\lambda)&=&\theta^3\left(\theta^{-2}-\lambda^{-2}\right)^2+
\left(\lambda^3-\theta^3\right)\left(\lambda^{-2}-\beta_\theta\right)^2+
\left(1-\lambda^3\right)(\beta_\theta-\beta_\lambda)^2\\&=&
\left(\frac{1}{\theta}-\frac{1}{\lambda}\right)\left(1-\frac{\theta}{\lambda}\right)
\left(1+\frac{\theta}{\lambda}\right)^2\\&&\mbox{}+
\left(\frac{1}{\theta}-\frac{1}{\lambda}\right)\left[\frac{\theta}{\lambda}+
\left(\frac{\theta}{\lambda}\right)^2+\left(\frac{\theta}{\lambda}\right)^3\right]
\left(1-\lambda^2\beta_\theta\right)^2\\&&\mbox{}+
\left(\frac{1}{\theta}-\frac{1}{\lambda}\right)\left(1-\frac{\theta}{\lambda}\right)
\theta\lambda^2\left(1-\lambda^3\right)
\left(\frac{\beta_\theta-\beta_\lambda}{\lambda-\theta}\right)^2.
\end{eqnarray*}
Since $8/9<1-\lambda^2\beta_\theta<1$ and, by (\ref{Es1}),
\[
0<\theta\lambda^2\left(1-\lambda^3\right)
\left(\frac{\beta_\theta-\beta_\lambda}{\lambda-\theta}\right)^2<\frac{1}{27},
\]
we conclude that
\[ 
G=\left(\theta^{-1}-\lambda^{-1}\right)^{-1}\!d_2^2(\theta,\lambda)=(1-z)(1+z)^2+
c_1(\theta,\lambda)\left(z+z^2+z^3\right)+c_2(\theta,\lambda)(1-z),
\]
with $z=\theta/\lambda\in(0,1)$, $8/9<c_1(\theta,\lambda)<1$ and
$0<c_2(\theta,\lambda)<1/27$. It follows that
\[
1<1+z-z^2-z^3+(8/9)\left(z+z^2+z^3\right)<G<1+2z+(1/27)(1-z)<3,
\]
which finally implies (\ref{Eq-Con2}). \cqfd\vspace{2mm}\\
It immediately follows from (\ref{Eq-Con1}) that the set $S_\eta=\{s_{\lambda_j},\;j\ge0\}$ with
$\lambda_j=(2j+1)2\eta^2/3$ is an $\eta$-net for the family $\overline{S}$ with respect to the Hellinger distance. On the other hand, given $\lambda\in(0,1/3)$ and $r\ge2\eta$, in order that
$s_{\lambda_j}\in{\cal B}(s_\lambda,r)$, it is required that $h^2(\lambda_j,\lambda)=
C(\lambda_j,\lambda)|\lambda_j-\lambda|< r^2$ which implies that
$|\lambda_j-\lambda|<(9/2)r^2$ and therefore
\[ 
|S_\eta\cap{\cal B}(s_\lambda,r)|\le1+(27/4)(r/\eta)^2\le
\exp\left[0.84(r/\eta)^2\right]\quad\mbox{for all }s_\lambda\in\overline{S}.
\]
It follows from Lemma~2 of Birg\'e (2006a) that $\overline{S}$ has a metric dimension
bounded by 3.4 and Corollary~3 of Birg\'e (2006a) implies that a suitable T-estimator
$\widetilde{s}$ built on $S_\eta$ has a risk satisfying 
\[
\sup_{0<\theta\le1/3}\Bbb{E}_{s_\theta}\left[h^2(s_\theta,\widetilde{s})\right]\le Cn^{-1}.
\]

Now setting $S_\eta=\{s_{\lambda_j},\;j\ge0\}$ with
$\lambda_j=\left(3+2j\eta^2/3\right)^{-1}$  we deduce as before that $S_\eta$ is an
$\eta$-net for $\overline{S}$ with respect to the $\Bbb{L}_2$-distance. In order that 
$s_{\lambda_j}\in{\cal B}(s_\lambda,x\eta)$, it is required that 
$d_2^2(\lambda_j,\lambda)=C(\theta,\lambda)|\lambda^{-1}_j-\lambda^{-1}|< x^2\eta^2$, 
which implies that $|\lambda^{-1}_j-\lambda^{-1}|<x^2\eta^2$. It follows that the number of 
elements of $S_\eta$ contained in the ball is bounded by $3x^2/2+1\le\exp\left(x^2/2\right)$ 
for $x\ge2$. Hence the metric dimension of $\overline{S}$ with respect to the 
$\Bbb{L}_2$-distance is bounded by $2$. It nevertheless follows from (\ref{Eq-A9}) that
the minimax risk over $\overline{S}$ is infinite when we use the $\Bbb{L}_2$-loss.

\subsection{Proof of Lemma~\ref{L-l2b}\labs{P1}}
Let us begin with a preliminary lemma.
%
\begin{lemma}\lab{L-l2a}
Let $F$ and $G$ be two disjoint sets with positive measures $\alpha=\mu(F)$ and
$\beta=\mu(G)$ and $g\in\overline{\Bbb{L}}_2$ such that $\inf_{x\in F}g(x)>0$. Set
$g_\varepsilon=g+\varepsilon(\alpha\1_{G} -\beta\1_{F})$ for $\varepsilon>0$.
Then $g_\varepsilon$ is a density for $\varepsilon$ small enough and for any $f\in
\overline{\Bbb{L}}_2$,
\begin{equation}
\lim_{\varepsilon\rightarrow0}\frac{1}{2\varepsilon}\left[d_2^2(g_\varepsilon,f)-
d_2^2(g,f)\right]=\alpha\int_G(g-f)\,d\mu-\beta\int_F(g-f)\,d\mu
\labe{Eq-L2a}
\end{equation}
and
\begin{equation}
\lim_{\varepsilon\rightarrow0}\frac{2}{\varepsilon}\left[h^2(g_\varepsilon,f)
-h^2(g,f)\right]=\beta\int_F\sqrt{fg^{-1}}\,d\mu-\alpha\int_G\sqrt{fg^{-1}}\,d\mu,
\labe{Eq-L2b}
\end{equation}
with the convention that $\int_G\sqrt{fg^{-1}}\,d\lambda=+\infty$ if either
$\mu(G\cap\{g=0\}\cap\{f>0\})>0$ or the integral diverges.
\end{lemma}
%
\noindent{\em Proof:}
Since $\int g_\varepsilon\,d\mu=1$ and $g_\varepsilon\ge0$ for $\varepsilon$
small enough $g_\varepsilon$ is a density. Moreover, setting
$k=\alpha\1_{G}-\beta\1_{F}$, we get
\[
d_2^2(g_\varepsilon,f)=\int(g+\varepsilon k-f)^2\,d\mu=
d_2^2(g,f)+\varepsilon^2\|k\|^2+2\varepsilon\int k(g-f)\,d\mu
\]
and (\ref{Eq-L2a}) follows. Let now $\Delta(\varepsilon)=
\varepsilon^{-1}\left[h^2(g_\varepsilon,f)-h^2(g,f)\right]$. Then
\begin{eqnarray*}
\Delta(\varepsilon)&=&\varepsilon^{-1}
\left[\int\sqrt{gf}\,d\mu-\int\sqrt{(g+\varepsilon k)f}\,d\mu\right]\\&=&
\varepsilon^{-1}\left[\int_F\left[\sqrt{gf}-\sqrt{(g-\varepsilon\beta)f}\right]d\mu
+\int_G\left[\sqrt{gf}-\sqrt{(g+\varepsilon\alpha)f}\right]d\mu\right]\\&=&
\int_F\frac{\beta\sqrt{f}}{\sqrt{g-\varepsilon\beta}+\sqrt{g}}\,d\mu-
\int_{G\cap\{g>0\}}\frac{\alpha\sqrt{f}}{\sqrt{g+\varepsilon\alpha}+\sqrt{g}}\,
d\mu\\&&\mbox{}-\int_{G\cap\{g=0\}\cap\{f>0\}}\sqrt{\alpha f/\varepsilon}
\,d\mu.
\end{eqnarray*}
When $\varepsilon$ tends to 0, the first integral converges to
$(\beta/2)\int_F\sqrt{fg^{-1}}\,d\mu$ and the second one converges to
$(\alpha/2)\int_{G\cap\{g>0\}}\sqrt{fg^{-1}}\,d\mu$, by monotone convergence.
The last one converges to $+\infty$ if $\mu(G\cap\{g=0\}\cap\{f>0\})>0$ and 0
otherwise, which achieves the proof of (\ref{Eq-L2b}).\cqfd\vspace{2mm}\\
If $\|v_1\vee v_2\|_\infty>2B$, we may assume, exchanging the roles of $v_1$ and $v_2$ if
necessary, that $\mu(A)>0$ with $A=\{v_1\ge v_2\;\;\mbox{and}\;\;v_1>2B\}$. Let
$C=\{v_1<B\wedge v_2\}$. If $\mu(C)>0$, we may apply Lemma~\ref{L-l2a} with
$F=A$, $G=C$, $g=v_1$ and $v'_1=g_\varepsilon$. We first set $f=t$. Since $v_1-t<B$ on
$C$ while $v_1-t>B$ on $A$, it follows from (\ref{Eq-L2a}) that $d_2(v'_1,t)<d_2(v_1,t)$ for
$\varepsilon$ small enough. If we now set $f=v_2$ and use (\ref{Eq-L2b}), we see that
$h(v'_1,v_2)<h(v_1,v_2)$ since $v_2\le v_1$ on $A$ and $v_2>v_1$ on $C$. We conclude by
setting $v'_2=v_2$. If $\mu(C)=0$, then $\mu(\{B\le v_1<v_2\})+\mu(\{v_2\le v_1<B\})=1$
and both sets have positive $\mu$-measure since $v_1\ne v_2$. In this case we set $F=\{B\le
v_1<v_2\}$, $G=\{v_2\le v_1\wedge u\}$ and $g=v_2$. Then $\mu(F)>0$ and $\mu(G)>0$
since $u\le B<v_2$ on $F$ and they are densities. If we use (\ref{Eq-L2a}) with $f=u$, we derive
that $d_2(v'_2,u)<d_2(v_2,u)$ for $\varepsilon$ small enough and if we use (\ref{Eq-L2b}) with
$f=v_1$, we derive that $h(v'_2,v_1)<h(v_2,v_1)$, in which case we set $v'_1=v_1$.

\subsection{Proof of Theorem~\ref{T-Selection}\labs{P3}}
We consider the family of tests $\psi(t_n,t_p,\bm{X})=\psi_{t_n,t_p,x}(\bm{X})$ provided by
the assumption with $x=A|p-n|$. Given this family of tests and $S=\{t_i, i\ge1\}$,
we define the random function ${\cal D}_{\bm{X}}$ on $S$ as in Birg\'e (2006a), i.e.\ we set 
${\cal R}_i=\{t_j\in S,j\ne i\,|\,\psi(t_i,t_j,\bm{X})  =t_j\}$ and
\begin{equation}
{\cal D}_{\bm{X}}(t_i)=\left\{\begin{array}{ll}{\displaystyle\sup_{t_j\in{\cal R}_i}
\left\{\st d(t_i,t_j)\right\}}&\!\!\mbox{if }\,{\cal R}_i\ne\emptyset;\\
\rule{0mm}{6mm}0&\!\!\mbox{if }\,{\cal R}_i=\emptyset.\end{array}\right.
\labe{Eq-H0}
\end{equation}
Given some $t_i\in S$, we want to bound
\[
\Bbb{P}_s\left[{\cal D}_{\bm{X}}(t_i)>xy_i\right]\quad\mbox{ for }x\ge1
\quad\mbox{and}\quad y_i=4d(s,t_i)\vee\sqrt{Aa^{-1}i2^i}.
\] 
Let us define the integer $K$ by $x^2<2^K\le2x^2$. Then 
\begin{equation}
K\ge1,\quad a2^{-i-K}(xy_i)^2\ge a2^{-i-1}y_i^2\ge Ai/2\quad\mbox{and}\quad
e^{-AK}\le x^{-2A/\log2}.
\labe{Eq-M32}
\end{equation}
Now, setting $y=xy_i$, observe that
\[
\Bbb{P}_s\left[{\cal D}_{\bm{X}}(t_i)>y\right]=
\Bbb{P}_s\left[\;\exists j\mbox{ with }d(t_i,t_j)>y
\mbox{ and }\psi(t_i,t_j,\bm{X})=t_j\right]\le\Sigma_1+\Sigma_2,
\]
with
\[
\Sigma_1=\sum_{j<i}\1_{d(t_i,t_j)>y}\;\Bbb{P}_s\left[\psi(t_i,t_j,\bm{X})=t_j\right];
\quad
\Sigma_2=\sum_{j>i}\1_{d(t_i,t_j)>y}\;\Bbb{P}_s\left[\psi(t_i,t_j,\bm{X})=t_j\right].
\]
If $i=1$, then $\Sigma_1=0$ and if $i\ge2$, we can use (\ref{Eq-tes2}) and $y\ge4d(s,t_i)$
to derive that
\begin{eqnarray*}
\Sigma_1&\le&B\sum_{j<i}\1_{d(t_i,t_j)>y}\;\exp\left[-a2^{-i}d^2(t_i,t_j)+A|i-j|
\right]\\&\le&B\exp\left[-a2^{-i}y_i^2x^2+Ai\right]\sum_{j\ge1}e^{-Aj}\\
&\le&B\frac{e^{-A}}{1-e^{-A}}\exp\left[-Ai\left(x^2-1\right)\right]\;\;\le\;\;B
\frac{e^{-A}}{1-e^{-A}}\exp\left[-A\left(x^2-1\right)\right]\\&\le&B
\left(1-e^{-A}\right)^{-1}\exp\left[-Ax^2\right]\;\;\le\;\;B
\left(1-e^{-A}\right)^{-1}x^{-2A/\log2},
\end{eqnarray*}
where we used (\ref{Eq-M32}), $i\ge1$ and $x\ge1$. Also, by (\ref{Eq-tes1}),
\begin{eqnarray*}
\Sigma_2&\le&B\sum_{j>i}\1_{d(t_i,t_j)>y}\;\exp\left[-a2^{-j}d^2(t_i,t_j)-A|i-j|\right]
\\&\le&B\sum_{j>i}\exp\left[-a2^{-j}y^2-A(j-i)\right]
\;\;=\;\;B\sum_{k=1}^{+\infty}\exp\left[-a2^{-i-k}y^2-Ak\right]
\\&\le&B\left[\sum_{k=1}^{K}\exp\left[-a2^{-i-k}y^2-Ak\right]+
\sum_{k>K}\exp[-Ak]\right]\;\;=\;\;B(\Sigma_3+\Sigma_4),
\end{eqnarray*}
with $\Sigma_4=e^{-AK}\left(e^A-1\right)^{-1}$ and, by  (\ref{Eq-M32}),
\begin{eqnarray*}
\Sigma_3&=&e^{-AK}\sum_{j=0}^{K-1}\exp\left[-a2^{-i-K+j}y^2+Aj\right]
\;\;\le\;\;e^{-AK}\sum_{j=0}^{K-1}\exp\left[-A(i2^{j-1}-j)\right]\\&\le&
e^{-AK}\sum_{j\ge0}\exp\left[-\left(2^{j-1}-j\right)\right]\;\;<\;\;3e^{-AK}.
\end{eqnarray*}
We finally get, putting all the bounds together and using (\ref{Eq-M32}) again,
\begin{equation}
\Bbb{P}_s\left[{\cal D}_{\bm{X}}(t_i)>xy_i\right]\le BC(A)x^{-2A/\log2}
\quad\mbox{for }x\ge1.
\labe{Eq-M31}
\end{equation}
As a consequence ${\cal D}_{\bm{X}}(t_i)<+\infty$ a.s.\ and we can define
\[
\widehat{s}_A=t_p\quad\mbox{with }p=\min\left\{j\,\left|\,{\cal D}_{\bm{X}}(t_j)
<\inf_i{\cal D}_{\bm{X}}(t_i)+\sqrt{Aa^{-1}}\right.\right\}.
\]
In view of the definition of ${\cal D}_{\bm{X}}$, $d(t_i,t_j)\le{\cal D}_{\bm{X}}(t_i)\vee
{\cal D}_{\bm{X}}(t_j)$, hence, for all $t_i\in S$, $d\left(\widehat{s}_A,t_i\right)\le
{\cal D}_{\bm{X}}(t_i)+\sqrt{Aa^{-1}}$ and
\[
d\left(\widehat{s}_A,s\right)\le{\cal D}_{\bm{X}}(t_i)+\sqrt{Aa^{-1}}+d(s,t_i)<
{\cal D}_{\bm{X}}(t_i)+y_i.
\]
It then follows from (\ref{Eq-M31}) that
\[
\Bbb{P}_s\left[d\left(\widehat{s}_A,s\right)>zy_i\right]\le BC(A)(z-1)^{-2A/\log2}
\quad\mbox{for }z\ge2.
\]
Integrating with respect to $z$ leads to
\[
\Bbb{E}_s\left[\left(d\left(\widehat{s}_A,s\right)/y_i\right)^q\right]\le BC(A,q)
\quad\mbox{for }1\le q<2A/\log2,
\]
and, since $t_i$ is arbitrary in $S$, 
\[
\Bbb{E}_s\left[d^q\left(\widehat{s}_A,s\right)\right]\le BC(A,q)\inf_{i\ge1}
\left[d^q(s,t_i)\vee\left(a^{-1}i2^i\right)^{q/2}\right]
\quad\mbox{for }1\le q<2A/\log2.\vspace{2mm}
\]
{\large{\bf Acknowledgements}}
I would like to thank Yannick Baraud for his remarks and the correction of some mistakes as well as two anonymous referees for their many suggestions.
\vspace{7mm}\\
\noindent{\large{\bf References\vspace{1mm}}}
\setlength{\parskip}{1mm}
{\small

BARAUD, Y. and BIRG\'E, L. (2011). Estimating composite functions by model selection. {\it
http://arxiv.org/abs/1102.2818}. To appear in {\it Ann. Inst. Henri Poincar\'e Probab.\ et Statist.} 

BIRG\'E, L. (1983). Approximation dans les espaces  m\'etriques   et th\'eorie de 
l'estimation. {\it Z. Wahrscheinlichkeitstheorie Verw. Geb.} {\bf65}, 181-237.  

BIRG\'E, L. (2004). Model selection for Gaussian regression with random design.
{\it Bernoulli} {\bf 10}, 1039 -1051.

BIRG\'E, L. (2006a). Model selection via testing : an alternative to (penalized) maximum
likelihood estimators. {\it Ann. Inst. Henri Poincar\'e Probab.\ et Statist.} {\bf 42}, 273-325.

BIRG\'E, L. (2006b). Statistical estimation with model selection. 
{\it Indagationes Mathematicae} {\bf 17}, 497-537.

BIRG\'E, L. (2007). Model selection for Poisson processes, in {\it Asymptotics: particles,
processes and inverse problems, Festschrift for Piet Groeneboom} (E.~Cator, G.~Jongbloed,
C.~Kraaikamp, R.~Lopuha\"a and J.~Wellner, eds). IMS Lecture Notes -- Monograph Series {\bf 55}, 32-64.

BIRG\'E, L. (2012). Robust tests for Model Selection. To appear in {\it From Probability to Statistics and Back: High-Dimensional Models and Processes} (M.~Banerjee, F.~Bunea, J.~Huang, V.~Koltchinskii and M.~Mathuis, eds). IMS Collections, Volume 9. 

BIRG\'E, L. and MASSART, P. (1997). From model selection to adaptive estimation. In {\it 
Festschrift for Lucien Le Cam: Research Papers in Probability and Statistics} (D. Pollard, E.
Torgersen and G. Yang, eds.), 55-87. Springer-Verlag, New York. 

BIRG\'E, L. and MASSART, P. (1998). Minimum contrast estimators on sieves: exponential 
bounds and rates of convergence. {\it Bernoulli} {\bf 4}, 329-375.

BIRG\'E, L. and MASSART, P.  (2000). An adaptive compression algorithm in Besov
spaces.  {\it Constructive Approximation} {\bf 16} 1-36.

BIRG\'E, L. and MASSART, P. (2001). Gaussian model selection. {\it J. Eur. Math. Soc.}
{\bf 3}, 203-268.

BIRG\'E, L. and ROZENHOLC, Y. (2006). How many bins should be put in a regular
histogram. {\it ESAIM-Probab. \& Statist.} {\bf 10}, 24-45.

CENCOV, N.N.  (1962). Evaluation of an unknown distribution density from observations. 
{\it Soviet Math.}  {\bf3}, 1559-1562.

DELYON, B. and JUDITSKY, A. (1996). On minimax wavelet estimators. {\it Appl.
Comput. Harmonic Anal.} {\bf 3}, 215-228.

DeVORE, R.A.\ and LORENTZ, G.G. (1993). {\it Constructive Approximation}. 
Springer-Verlag, Berlin. 

DEVROYE, L.  (1987). {\it A Course in Density Estimation}. Birkh\"auser, Boston.

DEVROYE, L.  and GY\"ORFI, L. (1985). {\it Nonparametric Density Estimation: The $L_1$
View}. John Wiley, New York.

DONOHO, D.L., JOHNSTONE, I.M., KERKYACHARIAN, G. and PICARD, D. (1996). 
Density estimation by wavelet thresholding. {\it Ann. Statist.} {\bf 24}, 508-539.

DONOHO, D.L. and LIU, R.C. (1987). Geometrizing rates of convergence I. Technical report 
137. Department of Statistics, University of California, Berkeley.  

EFROMOVICH, S. (2008). Adaptive estimation of and oracle inequalities for probability densities and characteristic functions.  {\it Ann. Statist.} {\bf 36}, 1127-1155.

JUDITSKY, A., RIGOLLET, P. and TSYBAKOV, A. (2008). Learning by mirror averaging.
{\it Ann. Statist.} {\bf 36}, 2183-2206.

KERKYACHARIAN, G. and PICARD, D. (2000). Thresholding algorithms, maxisets and
well-concentrated bases. {\it Test}   {\bf 9}, 283-344. 

Le CAM, L.M. (1973). Convergence of estimates under dimensionality restrictions. {\it
Ann. Statist.}  {\bf1} , 38-53. 

Le CAM, L.M. (1975). On local and global properties in the theory of asymptotic normality of
experiments. {\it Stochastic Processes and Related Topics, Vol.\ 1} (M. Puri, ed.), 13-54. 
Academic Press,  New York. 

Le CAM, L.M. (1986). {\it Asymptotic Methods in Statistical Decision Theory}. 
Springer-Verlag, New York. 

LOUNICI, K. (2008). Aggregation of density estimators for the $L^\pi$ risk with 
$1\le\pi<\infty$.  Unpublished manuscript.

MASSART, P.  (2007). Concentration Inequalities and Model Selection. In {\it  Lecture on
Probability Theory and Statistics, Ecole d'Et\'e de  Probabilit\'es de Saint-Flour XXXIII -
2003} (J.~Picard, ed.). Lecture Note in Mathematics, Springer-Verlag,  Berlin. 

PINSKER, M.S. (1980). Optimal filtration of square-integrable signals in Gaussian noise.
{\it  Problems  of Information Transmission} {\bf16}, 120-133.  

REYNAUD-BOURET, P., RIVOIRARD, V. and TULEAU-MALOT, C. (2011). Adaptive density estimation: A curse of support~? {\it J. Statist. Planning and Inference} {\bf 141}, 115-139.

RIGOLLET, T. (2006). Ph.D.\ thesis, University Pierre et Marie Curie, Paris.

RIGOLLET, T.\ and TSYBAKOV, A.B. (2007). Linear and convex aggregation of density
estimators. {\it Math. Methods of Statist.} {\bf 16}, 260-280.

SAMAROV, A. and TSYBAKOV, A.B. (2007) Aggregation of density estimators and dimension
reduction. {\it Advances in Statistical Modeling and Inference. Essays in Honor of Kjell A.
Doksum} (V.Nair, ed.), World Scientific, Singapore e.a., 233-251.

YANG, Y. (2000).  Mixing strategies for density estimation.  {\it Ann. Statist.}  {\bf28},
75-87. 

YANG, Y. and BARRON, A.R. (1998).  An asymptotic property of model selection criteria. 
{\it IEEE Transactions on Information Theory} {\bf 44}, 95-116.

YU, B. (1997). Assouad, Fano and Le~Cam. In {\it Festschrift for Lucien Le Cam: Research
Papers in Probability and Statistics} (D. Pollard, E. Torgersen and G. Yang, eds.), 423-435.
Springer-Verlag, New York. 
}
\vspace{8mm}\\
Lucien BIRG\'E\\ UMR 7599 ``Probabilit\'es et mod\`eles al\'eatoires"\\
Laboratoire de Probabilit\'es, bo\^{\i}te 188\\
Universit\'e Paris VI, 4 Place Jussieu\\
F-75252 Paris Cedex 05\\
France\vspace{2mm}\\
e-mail: lucien.birge@upmc.fr

\end{document}